\documentclass[a4paper]{amsart} 
\usepackage{amsmath}
\usepackage{amsthm}
\usepackage{amssymb} 
\usepackage{verbatim}
\usepackage[all]{xy}
\usepackage{multirow}
\usepackage{enumerate}

\usepackage{ytableau}


\newtheorem{l1}{Lemma}[section]
\newtheorem{c1}[l1]{Corollary}
\newtheorem{t1}[l1]{Theorem}
\newtheorem{p1}[l1]{Proposition}

\theoremstyle{definition}

\newtheorem{d1}[l1]{Definition}
\newtheorem{e1}[l1]{Example}

\theoremstyle{remark}
\newtheorem{r1}[l1]{Remark}

\renewcommand{\mod}{\operatorname{\mathsf{mod}}\nolimits}

\newcommand{\ind}{\operatorname{\mathsf{ind}}\nolimits}
\newcommand{\Coker}{\operatorname{\mathsf{Coker}}\nolimits}
\newcommand{\Ker}{\operatorname{\mathsf{Ker}}\nolimits}
\newcommand{\Ext}{\operatorname{\mathsf{Ext}}\nolimits}
\newcommand{\End}{\operatorname{\mathsf{End}}\nolimits}
\renewcommand{\Im}{\operatorname{\mathsf{Im}}\nolimits}
\newcommand{\sub}{\operatorname{\mathsf{sub}}\nolimits}

\newcommand{\soc}{\operatorname{soc}\nolimits}
\renewcommand{\top}{\operatorname{top}\nolimits}

\newcommand{\add}{\operatorname{\mathsf{add}}\nolimits}

\title{The Lattice of Subobject Closed Subcategories and Colocal Type}

\author{Apolonia Gottwald}

\address{Apolonia Gottwald, Max Planck Institute for Mathematics, Vivatgasse 7, D-53111 Bonn, Germany}

\date{\today}

\email{alogisma@gmx.de}

\begin{document}
	
\begin{abstract}
	In this paper, we consider abelian length categories, a generalization of module categories over Artin algebras. Let $\mathcal{A}$ be an abelian length category of colocal type. We show that the lattice $\mathsf{S}(\mathcal{A})$ of full additive subobject closed subcategories of $\mathcal{A}$ is distributive. Furthermore, we give a characterization of abelian length categories of colocal type. 
	
	If $A$ is an algebra of colocal type over an algebraically closed field, then this characterization is especially simple and we can describe the lattice $\mathsf{S}(\mod A)$ up to isomorphism.
\end{abstract}
	
\maketitle

\tableofcontents

\section{Introduction}

Let $\mathcal{A}$ be abelian length category (a generalization of module categories over Artin algebras). In this paper, we give a characterization of abelian length categories of colocal type. Furthermore, we are interested in the lattice $\mathsf{S}(\mathcal{A})$ of full additive subobject closed subcategories of $\mathcal{A}$. In particular, we show that $\mathsf{S}(\mathcal{A})$ is distributive if $\mathcal{A}$ is of colocal type.

Subobject closed subcategories have not yet been extensively studied, but there are connections to different parts of representation theory. For example, let $A$ be a finite dimensional algebra: then every infinite submodule closed subcategory of $\mod A$ contains a minimal infinite submodule closed category, as Ringel proved in \cite{Rin}.

Krause and Prest have used submodule closed subcategories in \cite{KP} to show that there is a filtration of the Ziegler spectrum that is indexed by the Gabriel-Roiter filtration. 

Furthermore, if $A$ is a hereditary Artin algebra, then there is a natural bijection between the elements of the Weyl group of $A$ and the full, additive cofinite submodule closed subcategories of $\mod A$. This has been proved by Oppermann, Reiten and Thomas in \cite{ORT} for algebras over finite and algebraically closed fields and in general the author's preprint \cite{Got}. 

In Section \ref{SecColoc}, we see that $\mathsf{S}(\mathcal{A})$ is distributive if $\mathcal{A}$ is of colocal type. Algebras of colocal type have been studied repeatedly: for example, a first characterization dates back to Tachikawa in 1959, see \cite{Tac}; two gaps in the proof were filled by Sumioka in 1984, see \cite{Sum}.

In this paper, we give a new characterisation for abelian length categories of colocal type. Note that we are equating objects with isomorphism classes of objects; in particular, all sums over simple objects are actually sums over isomorphism classes of simple objects.

We define

\begin{d1}
	For all simple objects $S, T \in \mathcal{A}$ let
	\begin{equation*}
	d^1_S(S, T) := \dim_{\End(S)^{op}} \Ext^1(S, T)
	\end{equation*}
	and
	\begin{equation*}
	d^1_T(S, T) := \dim_{\End(T)} \Ext^1(S, T).
	\end{equation*}	
\end{d1}

Then we can show the following:

\begin{t1} \label{ThShape}
	An abelian length category $\mathcal{A}$ is of colocal type if and only if the following conditions hold for all simple objects $S \in \mathcal{A}$:

		\begin{equation*} \tag{C1}
		\sum_{T \text{ simple}} d^1_T(S, T) \le 1
		\end{equation*}
		
		\begin{equation*} \tag{C2}
		\sum_{T \text { simple}} d^1_T(T, S) \le 2
		\end{equation*}
	\begin{enumerate}[\textup{(C3)}]
		\item If there is a simple object $S'$ with $\Ext^1(S,S') \neq 0$, let 	\begin{equation*}
		\mathcal{T} := \{T \text{ simple and }\Ext^1(T, S) \neq 0  \mid \exists Z: l(Z) = 3, \soc Z = S',  \top Z = T\}
		\end{equation*}
		 Then
		\begin{equation*}
		\sum_{T \in \mathcal{T}} d^1_T(T, S) \le 1.
		\end{equation*}
	\end{enumerate}
\end{t1}
While the last condition is more complicated then the first two, there are several ways to state it, see Proposition \ref{PropEquiv}. In particular, it is often equivalent to a condition on the 2-extensions between simple objects.

Furthermore, for algebras of colocal type over algebraically closed fields, we can completely describe the lattice $\mathsf{S}(\mod A)$. 

This paper is organised in the following way: In Section \ref{SecIndecObj}, we show that the distributivity of $\mathsf{S}(\mathcal{A})$ is equivalent to a condition on the submodule relations in $\mathcal{A}$.
We can show in the next section that the following is an even stronger property: every subobject of an indecomposable object in $\mathcal{A}$ is itself indecomposable. Such categories are said to be of colocal type.

We characterize these categories in Section \ref{SecExt} to \ref{SecEquiv}: First, we show that two conditions on the $\Ext$-quiver hold. Weaker conditions hold if $\mathsf{S}(\mathcal{A})$ is distributive. In Section \ref{SecThird}, we give different formulations and a proof of the third condition that abelian length categories of colocal type fulfil. Again, we see that a weaker condition is fulfilled if $\mathsf{S}(\mathcal{A})$ is distributive. 

In Section \ref{SecEquiv}, we prove that every abelian length category which fulfils the three conditions is of colocal type.

Returning to the lattice $\mathsf{S}(\mathcal{A})$, we show in the next section that it is the Cartesian product of certain sublattices.

In Section \ref{SecString}, we give a short summary of some facts about string algebras and their modules, which we need in the next and final section.

There, we assume that $\mathcal{A} \equiv \mod kQ/I$ for some field $k$ and some quiver $Q$ with an admissible ideal $I$. In this case, $A$ is of colocal type if and only if $A$ is a string algebra and no vertex in $Q$ is starting point of more than one arrow. For these algebras, we get a complete, explicit description of the lattice $\mathsf{S(mod} ~A)$.

Throughout this paper, we use the notation $X \mid Y$ if the object $X$ is a direct summand of $Y$ and $X \nmid Y$ if $X$ is not a direct summand of $Y$. For the length of $X$, we write $l(X)$.

Furthermore, we use a matrix notation
\begin{equation*}
\begin{bmatrix}
f_{11} & \dots & f_{1m}\\
&\dots \\
f_{n1} & \dots & f_{nm}
\end{bmatrix}: X_1 \oplus \dots \oplus X_{m} \rightarrow Y_1 \oplus \dots \oplus Y_n
\end{equation*}
to denote morphisms $f_{ij}: X_i \rightarrow Y_j$ for $1 \le i \le m$ and $1 \le j \le n$.

\section{Conditions on indecomposable objects} \label{SecIndecObj}
In this section we characterize the abelian length categories $\mathcal{A}$ with distributive lattices $\mathsf{S}(\mathcal{A})$ in terms of the subobject relations between the objects of $\mathcal{A}$.

We start with the definition of a distributive lattice, as given for example in \cite{Ste}, p. 69:

\begin{d1} \label{DefDistributive}
	A lattice $L$ is called \textit{distributive} if 
	\begin{equation*}
	(a \vee b) \wedge c = (a \wedge c) \vee (b \wedge c)
	\end{equation*}
	for all $a, b, c \in L$.
\end{d1}

Now let $\mathsf{S}(\mathcal{A})$ be the set of full additive subobject closed subcategories as in \cite{KP}. It is partially ordered by inclusion and a complete lattice. 

The join $a \vee b$ for two categories $a, b \in \mathsf{S}(\mathcal{A})$ is the smallest subcategory in $\mathsf{S}(\mathcal{A})$ which contains both $a$ and $b$. The meet $a \wedge b$ is the largest category in $\mathsf{S}(\mathcal{A})$ that is contained in both $a$ and $b$. 

The meet coincides with the intersection $a \cap b$: all subobjects of direct sums of objects in $a \cap b$ are again objects in $a \cap b$, since $a$ and $b$ are subobject closed. The join consists of all subobjects of direct sums of objects in $a$ and $b$.

Every category in $\mathsf{S}(\mathcal{A})$ is completely determined by the isomorphism classes of indecomposable objects it contains.

For a class $\mathcal{X}$ of objects let $\sub \mathcal{X}$ be the category that consists of all subobjects of direct sums of objects in $\mathcal{X}$. This is the smallest category in $\mathsf{S}(\mathcal{A})$ that contains $\mathcal{X}$. Furthermore, let $\sub X := \sub \{X\}$.\\

In the following case, $\mathsf{S}(\mathcal{A})$ is not distributive:

\begin{l1} \label{LemDistributive}
	If there exists an indecomposable object $X \in \mathcal{A}$, and objects $Y_1, Y_2 \in \mathcal{A}$ so that $X \in \sub Y_1 \vee \sub Y_2$ but $X \notin \sub Y_i$ for all $1 \le i \le 2$, then
	\begin{equation*}
	(\sub Y_1 \vee \sub Y_2) \wedge \sub X \neq (\sub Y_1 \wedge \sub X) \vee (\sub Y_2 \wedge \sub X).
	\end{equation*}
\end{l1}

\begin{proof}
	By the assumption 
	\begin{equation*}
	X \in  \sub Y_1 \vee \sub Y_2
	\end{equation*}
	and by definition
	\begin{equation*}
	X \in \sub X,
	\end{equation*}
	so 
	\begin{equation*}
	X \in \left( \sub Y_1 \vee \sub Y_2 \right) \wedge \sub X.
	\end{equation*}
	But
	\begin{equation*}
	X \notin (\sub Y_1 \wedge \sub X) \vee (\sub Y_2 \wedge \sub X),
	\end{equation*}
	since otherwise there were some objects 
	\begin{equation*}
	X_i \in \sub Y_i \wedge \sub X
	\end{equation*}
	for $1 \le i \le 2$ with a monomorphism
	\begin{equation*}
	f: X \rightarrowtail  X_1 \oplus X_2.
	\end{equation*}
	We can assume without loss of generality that the components $X \rightarrow X_1$ and $X \rightarrow X_2$ of $f$ are epimorphisms. Since $X \notin \sub Y_i$ for $1 \le i \le 2$, these are not isomorphisms and $l(X_i) < l(X)$. 
	
	Since $X_1 \oplus X_2 \in \sub X$, there is some $\alpha \in \mathbb{N}$ with a monomorphism
	\begin{equation*}
	X_1 \oplus X_2 \rightarrowtail X^{\alpha} \rightarrowtail X_1^{\alpha} \oplus X_2^{\alpha}.
	\end{equation*}
	The restriction of this concatenation to $X_1$ is a monomorphism
	\begin{equation*}
	\begin{bmatrix}
	g_1 \\g_2
	\end{bmatrix}: X_1 \rightarrowtail X_1^{\alpha} \oplus X_2^{\alpha}.
	\end{equation*}
	
	Then $g_1$ cannot be a monomorphism, since this would imply $X_1 \mid X^{\alpha}$ and thus $X_1 \cong X^{\beta}$ for some $\beta \le \alpha$, since $X$ is indecomposable. But this is a contradiction to $X \notin \sub X_1$. So $l(\Im(g_1)) < l(X_1)$. 

	There is a concatenation of monomorphisms
	\begin{equation*}
	X \rightarrowtail X_1 \oplus X_2 \rightarrowtail \Im(g_1) \oplus \Im(g_2) \oplus X_2 \rightarrowtail \Im(g_1) \oplus X^{\alpha+1}.
	\end{equation*}
	With $X \notin \sub X_2$, we have $0 \neq \Im(g_1)$. We can set $X'_1 := \Im(g_1)$ and $X'_2 := X^{\alpha+1}$. Then there is a monomorphism $X \rightarrowtail X'_1 \oplus X'_2$ with $X \notin \sub X'_1$ and $X \notin X'_2$. Since this is the same situation as before, we get inductively an infinite sequences of non-zero objects $X_1, X'_1, X^{(2)}_1, X^{(3)}_1\dots $ with
	\begin{equation*}
	l(X_1) > l(X'_1) > X^{(2)}_1 > X^{(3)}_1 > \dots > 0.
	\end{equation*} 
	This cannot be true, since $l(X_1)$ is finite.
	So one of these objects is either the zero object or a direct sum of copies of $X$. Thus $X \in \sub X_1$ or $X \in \sub X_2$ must hold 
	and the proof is complete.
\end{proof}

In fact, we get the following equivalence:

\begin{p1} \label{PropDistributive}
	The following statements are equivalent:
	\begin{enumerate}[\textup{(}1\textup{)}]
		\item The lattice $\mathsf{S}(\mathcal{A})$ is distributive
		
		\item If $X \in \mathcal{A}$ is indecomposable and there are objects $Y_1, Y_2 \in \mathcal{A}$, so that $X \in  \sub Y_1 \vee \sub Y_2$ then $X \in \sub Y_i$ for some $1 \le i \le 2$.
		
		\item For all index sets $I$ and categories $a_i \in \mathsf{S}(\mathcal{A})$, $i \in I$ we have 
		\begin{equation*}
		\ind (\bigvee_{i\in I} a_i) = \bigcup_{i \in I} \ind a_i.
		\end{equation*}
		
		\item For all $a, b \in \mathsf{S}(\mathcal{A})$, we have
		\begin{equation*}
		\ind(a \vee b) = \ind a \cup \ind b.
		\end{equation*}
	\end{enumerate}
\end{p1}

\begin{proof}
	$(1) \Rightarrow (2)$ is the result of Lemma \ref{LemDistributive}.
	
	$(2) \Rightarrow (3)$: The direction 
	\begin{equation*}
	\ind (\bigvee_{i\in I} a_i) \supseteq \bigcup_{i \in I} \ind a_i
	\end{equation*}
	is clear.
	For the other direction, we look at an object
	\begin{equation*}
	X \in \ind \left(\bigvee_{i \in I} a_i\right).
	\end{equation*}
	There are objects $A_i \in a_i$ with a monomorphism
	\begin{equation*}
	X \rightarrowtail \bigoplus_{i \in I} A_i
	\end{equation*}
	and thus
	\begin{equation*}
	X \in \bigvee_{i \in I} \sub A_i.
	\end{equation*}
	The object $\bigoplus_{i \in I} A_i$ must be of finite length; thus $A_i = 0$ for all except finitely many $i \in I$. With (2) and an induction, we get  
	\begin{equation*}
	X \in \sub A_i
	\end{equation*}
	for at least one ${i \in I}$ and thus
	\begin{equation*}
	X \in a_i.
	\end{equation*}
	So
	\begin{equation*}
	X \in \bigcup_{i \in I} \ind a_i
	\end{equation*}
	and
	\begin{equation*}
	\ind (\bigvee_{i \in I} a_i) = \bigcup_{i \in I} \ind a_i.
	\end{equation*}
	
	$(3) \Rightarrow (4)$ is clear.
	
	$(4) \Rightarrow (1)$: Let $a, b, c \in \mathsf{S}(\mathcal{A})$. Then
	\begin{align*}
	\ind ((a \vee b) \wedge c) 
	&= (\ind a \cup \ind b) \cap \ind c \\
	&= (\ind a \cap \ind c) \cup (\ind a \cap \ind c)\\
	&= \ind((a \wedge c) \vee (a\wedge c)).
	\end{align*}
	Since $a, b, c$ are completely determined by their indecomposable objects, 
	\begin{equation*}
	(a \vee b) \wedge c = (a \wedge b) \vee (a \wedge c)
	\end{equation*}
	and $\mathsf{S}(\mathcal{A})$ is distributive.
\end{proof}

We can generalize the notion of a distributive lattice as e.g. in \cite{IK}, p. 1227: 

\begin{d1}
	A complete lattice $\Lambda$ is a \textit{frame} if for all index sets $I$ and elements $a, b_i$ with $i \in I$ the equation
	\begin{equation*}
	a \wedge \left( \bigvee_{i \in I} b_i \right) = \bigvee_{i \in I} (a \wedge b_i)
	\end{equation*}
	holds.
\end{d1}

Obviously, every frame is also distributive. But in general, not every distributive lattice is a frame. An exception are lattices of subobject closed categories: 

\begin{c1}
	The lattice $\mathsf{S}(\mathcal{A})$ is distributive if and only if it is a frame.
\end{c1}

\begin{proof}
	This follows from part (3) of Proposition \ref{PropDistributive}.
\end{proof}

\section{Categories of colocal type} \label{SecColoc}

\begin{d1}
	Let $\mathcal{A}$ be a finite length category. An object $X \in \mathcal{A}$ is called \textit{colocal}, if its socle is simple. Equivalently, it is colocal if every non-zero subobject of $X$ is indecomposable.
	
	Dually, $X$ is \textit{local}, if its top is simple or equivalently, if every non-zero factor object of $X$ is indecomposable.
\end{d1}

\begin{d1}
	We call a category $\mathcal{A}$ \textit{of colocal type} if every indecomposable object in $\mathcal{A}$ is colocal. If there is some Artin algebra $A$ so that $\mathcal{A} = \mod A$, then we also say that $A$ is of colocal type.
\end{d1}

For these categories, $\mathsf{S}(\mathcal{A})$ is always distributive. To show this, we need the following lemma:

\begin{l1} \label{LemTwoMonos}
\begin{enumerate}[(a)]
	\item 	If there are objects $V_1, V_2, X$ with a monomorphism
	\begin{equation*}
\begin{bmatrix}
f_1 & f_2
\end{bmatrix}: V_1 \oplus V_2 \rightarrowtail X,
\end{equation*}
then there is also a monomorphism
\begin{equation*}
X \rightarrowtail \Coker f_1 \oplus \Coker f_2.
\end{equation*}

\item If there are objects $X, Y_1, Y_2$ with a monomorphism
	\begin{equation*}
\begin{bmatrix}
f_1 \\ f_2
\end{bmatrix}: X \rightarrowtail Y_1 \oplus Y_2,
\end{equation*}
then there is also a monomorphism
\begin{equation*}
\Ker f_1 \oplus \Ker f_2 \rightarrowtail X.
\end{equation*}
\end{enumerate}
\end{l1}

\begin{proof}
	First, we prove that (a) holds: Under the assumptions on $V_1, V_2, X$, there is an exact sequence
	\begin{equation*}
	\xymatrix {	
		0 \ar[r]&  V_1 \oplus V_2 \ar^{\left[\begin{smallmatrix}
			f_1 & 0\\0 & f_2
			\end{smallmatrix}\right]}[r] & X \oplus X \ar[r]^-{g}&\Coker f_1 \oplus \Coker f_2 \ar[r]&0
	}
	\end{equation*} 
	for some morphism $g$ with $\Ker g \cong V_1 \oplus V_2$. Since there exists a monomorphism $\Ker g \rightarrowtail X$,	$g$ induces a monomorphism
	\begin{equation*}
	X \rightarrowtail \Coker f_1 \oplus \Coker f_2.
	\end{equation*}
The proof of (b) is similar:
	Under these assumptions, there is a monomorphism $f: X \rightarrowtail \Im(f_1) \oplus \Im(f_2)$ with an exact diagram
\begin{equation*}
\xymatrix@C=3pc {
	&0 \ar[d]&0\ar[d]&0\ar[d]\\
	0 \ar[r]&  0 \ar[d]\ar[r] & X \ar@{=}[r] \ar[d]& X \ar[r]\ar^{f}[d]&0\\
	0 \ar[r]&  \Ker f_1 \oplus \Ker f_2 \ar@{=}[d]\ar[r] & X \oplus X \ar^-{\left[\begin{smallmatrix}
		f_1 & 0\\0 & f_2
		\end{smallmatrix}\right]}[r] \ar[d]&\Im(f_1) \oplus \Im(f_2) \ar[r]\ar[d]&0\\
	& \Ker f_1 \oplus \Ker f_2 \ar[d] & X \ar[d]&\Coker f' \ar[d]&\\								
	&0&0&0			
}.
\end{equation*}
So there is an exact sequence
\begin{equation*}
\xymatrix{0 \ar[r]& \Ker f_1 \oplus \Ker f_2 \ar[r] & X \ar[r]&\Coker f \ar[r]& 0}
\end{equation*}
 and in particular, there is a monomorphism $\Ker f_1 \oplus \Ker f_2 \rightarrowtail X$.
\end{proof}

Now we can prove that $\mathsf{S}(\mathcal{A})$ is distributive if $\mathcal{A}$ is of colocal type:

\begin{p1} \label{LemOneDirection}
	If  $\mathcal{A}$ is of colocal type, then $\mathsf{S}(\mathcal{A})$ is distributive. Furthermore, if $X$ is indecomposable, then for all $Y_1, Y_2 \in \mathcal{A}$ with a monomorphism $X \rightarrowtail Y_1 \oplus Y_2$, there is a monomorphism $X \rightarrowtail Y_1$ or $X \rightarrowtail Y_2$.
\end{p1}

\begin{proof}
	If  $\mathcal{A}$ is of colocal type, then every subobject of an indecomposable object is indecomposable.
	Let $X, Y_1, Y_2$ be as in the assumptions. Then there is a monomorphism 
	\begin{equation*}
	\begin{bmatrix}
	f_1 \\ f_2
	\end{bmatrix}: X \rightarrowtail Y_1 \oplus Y_2
	\end{equation*}
	and by Lemma \ref{LemTwoMonos}, there is a monomorphism $\Ker f_1 \oplus \Ker f_2 \rightarrowtail X$.

	So $\Ker f_1 = 0$ or $\Ker f_2 = 0$ and $f_1$ or $f_2$ is a monomorphism. By Proposition \ref{PropDistributive}, the lattice $\mathsf{S}(\mathcal{A})$ is distributive.
\end{proof}

On the other hand, we have the following:

\begin{l1} \label{LemSoc}
	If $\mathsf{S}(\mathcal{A})$ is distributive, then the socle of any indecomposable object is of the form $S^m$ for some simple object $S$ and some $m \in \mathbb{N}$.
\end{l1}

\begin{proof}
	Let $X \in \mathcal{A}$ be indecomposable with simple objects $S_1 \ncong S_2$ so that $S_1 \oplus S_2 \mid \soc X$. Then there are monomorphisms $f_i: S_i \rightarrowtail X$ and by Lemma \ref{LemTwoMonos}, there is a monomorphism $X \rightarrowtail \Coker f_1 \oplus \Coker f_2$.
	
	We prove that  $\mathsf{S}(\mathcal{A})$ is not distributive with induction on $l(X)$.
	
	First, suppose that $l(X) = 3$. Then every indecomposable direct summand $Y$ of $\Coker f_1 \oplus \Coker f_2$ has length at most 2. So the socle of $Y$ is simple and either $S_1$ or $S_2$ is not a subobject of $Y$. Thus, $X \notin \sub Y$. By Proposition \ref{PropDistributive}, $\mathsf{S}(\mathcal{A})$ is not distributive.
	
	Now let the assumption be proved for all indecomposable objects with length smaller that $X$. If there is some indecomposable direct summand $Y$ of $\Coker f_1 \oplus \Coker f_2$ so that $X \in \sub Y$, then $T_1 \oplus T_2 \mid \soc Y$. Since $l(Y) < l(X)$, the lattice $\mathsf{S}(\mathcal{A})$ is not distributive by the inductive assumption.
\end{proof}

To prove our main result about colocal abelian length categories, we need the following lemma:

\begin{l1} \label{LemUniserial}
	
	If $\mathcal{A}$ is not of colocal type, then there are objects $V_1, V_2$, non-simple, colocal objects $Y_1, Y_2, \dots, Y_m$, an indecomposable object $X$ and a simple object $S$ with exact sequences 
	\begin{equation} \label{EqUniserial1}
	\xymatrix{0 \ar[r] & V_1 \oplus V_2 \ar^-{\left[\begin{smallmatrix}
	f_1 \\ f_2
	\end{smallmatrix}\right]}[r] & X \ar[r] & S \ar[r] & 0}
	\end{equation}
	and
	\begin{equation} \label{EqUniserial2}
	\xymatrix{0 \ar[r] & X \ar[r] & \bigoplus_{i=1}^m Y_i \ar[r] & S \ar[r] & 0}.
	\end{equation}
	For such objects, the following sequences are exact	for $1 \ le i, j \in 2$ and $i \neq j$:
	\begin{equation} \label{EqUniserial3}
	\xymatrix{0 \ar[r] & V_j \ar[r] & \Coker f_i \ar[r] & S \ar[r]& 0}.
	\end{equation}

\end{l1}

\begin{proof}
	If $\mathcal{A}$ is not of colocal type, then there is an indecomposable object $X$ that is not colocal. So there are objects $V_1 \neq 0 \neq V_2$ with a monomorphism 
	\begin{equation*}
	f= \begin{bmatrix}
	f_1 & f_2
	\end{bmatrix}: V_1 \oplus V_2 \rightarrowtail X.
	\end{equation*}
	Let $S$ be a simple factor module of $\Coker f$. Then there is some $V$ with $V_1 \oplus V_2 \subseteq V$ and an exact sequence 
	\begin{equation*}
	\xymatrix{0 \ar[r] & V \ar[r] & X \ar[r] & S \ar[r] & 0}.
	\end{equation*}
	If $V$ is indecomposable, then it is of smaller length than $X$ and not colocal. 
	
	Inductively, we can assume that $\Coker f = S$ and $V_1 \oplus V_2 = V$. By Lemma \ref{LemTwoMonos}, there is a monomorphism
	\begin{equation*}
	g: X \rightarrowtail \Coker f_1 \oplus \Coker f_2.
	\end{equation*}

	The following diagram is exact for all $1 \le i, j \le 2$, $i \neq j$, since all columns and the first and second row are exact:
	\begin{equation*}
	\xymatrix{
		& 0 \ar[d] & 0 \ar[d] & 0 \ar[d]\\
		0 \ar[r] & V_i \ar@{=}[r] \ar[d] & V_i \ar[r] \ar[d] & 0 \ar[r] \ar[d] & 0\\
		0 \ar[r] & V_1 \oplus V_2 \ar[r] \ar[d] & X \ar[r] \ar[d]& S \ar[r] \ar@{=}[d]& 0\\
		0 \ar[r]& V_j \ar[r] \ar[d] & \Coker f_i \ar[r] \ar[d] & S \ar[r] \ar[d] & 0\\
		& 0& 0& 0&
	}.
	\end{equation*}
	Thus, the following diagram is exact, since all columns and the first and second row are exact:
	\begin{equation*}
	\xymatrix{
		& 0 \ar[d] & 0 \ar[d] & 0 \ar[d]\\
		0 \ar[r] & V_1 \oplus V_2 \ar[r] \ar@{=}[d] &X \ar[r] \ar^g[d] & S \ar[r] \ar[d] &  0\\
		0 \ar[r] & V_1 \oplus V_2 \ar[r] \ar[d] & \Coker f_1 \oplus \Coker f_2 \ar[r] \ar[d]& S^2 \ar[r] \ar[d]& 0\\
		0 \ar[r]& 0 \ar[r] \ar[d] &  \Coker g\ar[r] \ar[d] & S \ar[r] \ar[d] & 0\\			
		& 0& 0& 0&
	}.
	\end{equation*}
	So $\Coker g = S$. Let
	\begin{equation*}
	\bigoplus_{i=1}^m Y_i := \Coker f_1 \oplus \Coker f_2
	\end{equation*}
	be a decomposition of $\Coker f_1 \oplus \Coker f_2$ into indecomposable direct summands. Then we get the exact sequences (\ref{EqUniserial1}) - (\ref{EqUniserial3}). Since $l(Y_i) < l(X)$, we can inductively find some object $X$ so that $Y_i$ is colocal for all $1 \le i \le m$.
\end{proof}

\section{Conditions on the $\Ext$-quiver} \label{SecExt}

For simple objects $S, T \in \mathcal{A}$, we define
\begin{equation*}
d^1_S(S, T) := \dim_{\End(S)^{op}} \Ext^1(S, T)
\end{equation*}
and
\begin{equation*}
d^1_T(S, T) := \dim_{\End(T)} \Ext^1(S, T).
\end{equation*}

In this section, we show that for an abelian length category $\mathcal{A}$ of colocal type all simple objects $S \in \mathcal{A}$ fulfil the following conditions:

		\begin{equation*} \tag{C1}
		\sum_{T \text{ simple}} d^1_T(S, T) \le 1
		\end{equation*}
		
		\begin{equation*} \tag{C2}
		\sum_{T \text { simple}} d^1_T(T, S) \le 2.
		\end{equation*}

 Weaker conditions hold if $\mathsf{S}(\mathcal{A})$ is distributive.

Recall that we can interpret $\Ext^n(S, T)$ as the group of equivalence classes of $n$-fold extensions of $S$ by $T$, see for example \cite{Ben}, Definition 2.6.1. For $\Ext^1(S, T)$, the abelian group structure corresponds to the Baer sum, see \cite{Zim}, Section 1.8.2:
	
	For
	\begin{equation*}
	\eta_1: \xymatrix{0 \ar[r] & X' \ar^{f_1}[r] & X_0 \ar^{g_1}[r] &X\ar[r] &0 }
	\end{equation*}
	and 
	\begin{equation*}
	\eta_2: \xymatrix{0 \ar[r] & X' \ar^{f_2}[r] & X'_0 \ar^{g_2}[r] &X\ar[r] &0},
	\end{equation*}
	there is an object $Z_1$ with a commutative diagram
	\begin{equation*}
	\xymatrix@C=3pc{0 \ar[r] & X' \oplus X'\ar[r] \ar@{=}[d]& Z_1 \ar[r] \ar[d]&X\ar[r]\ar^{\left[\begin{smallmatrix}
			\text{id}\\\text{id}
			\end{smallmatrix}\right]}[d] &0 \\
		0 \ar[r] & X' \oplus X' \ar^{\left[\begin{smallmatrix}
			f_1 & 0\\0 & f_2
			\end{smallmatrix}\right]}[r] & X_0 \oplus X'_0\ar^{\left[\begin{smallmatrix}
			g_1 & 0\\0 & g_2
			\end{smallmatrix}\right]}[r] &X \oplus X\ar[r] &0 }
	\end{equation*}
	and an object $Z$ with a commutative diagram
	\begin{equation*}
	\xymatrix{0 \ar[r] & X' \oplus X'\ar[r] \ar^{\left[\begin{smallmatrix}
			\text{id}&\text{id}
			\end{smallmatrix}\right]}[d]& Z_1 \ar[r] \ar[d]&X\ar[r]\ar@{=}[d] &0 \\
		0 \ar[r] & X' \ar[r] & Z\ar[r] &X\ar[r] &0 }.
	\end{equation*}
	The object $Z_1$ can be found by taking the pullback, while $Z$ can be found via the pushout.
	
	Then
	\begin{equation*}
	\xymatrix{	0 \ar[r] & X' \ar[r] & Z\ar[r] &X\ar[r] &0 }
	\end{equation*}
	is the exact sequence $\eta_1 + \eta_2$.\\
	
	Furthermore, we need the following result by Gabriel, see \cite{Gab}, p.81:
	
	\begin{t1}
		An abelian length category $\mathcal{A}$ is equivalent to the module category of an Artinian ring if and only if
		\begin{enumerate}
			\item $\mathcal{A}$ has only finitely many simple objects.
			
			\item  $d^1_T(S, T) < \infty$ for all simple objects $S, T \in \mathcal{A}$.
			
			\item The supremum of the Loewy lengths of the objects in $\mathcal{A}$ is finite.
		\end{enumerate}
	\end{t1}

In particular, if condition (2) is fulfilled, then for every finite set of simple objects $S_1, S_2, \dots, S_n$ and $r \in \mathbb{N}$, the subcategory of $\mathcal{A}$ that consists of the objects of Loewy length smaller or equal to $r$ with composition factors in $\{S_1, \dots, S_n\}$ is equivalent to such a module category.

If $\mathcal{A}$ contains simple objects $S, T$ with $\dim_{\End(T)}\Ext^1(S, T) = \infty$, then there is a subcategory $\mathcal{A}'_m$ of $\mathcal{A}$ so that $\dim_{\End(T)}\Ext^1_{\mathcal{A}'_m}(S, T) = m$ for all $m \in \mathbb{N}$.

\begin{l1} \label{LemTwoSimples}
	Let $\mathcal{A}$ be an abelian length category.
	\begin{enumerate}[(a)]
		\item If $\mathcal{A}$ is of colocal type, then (C1) is fulfilled.
		
		\item If $\mathsf{S}(\mathcal{A})$ is distributive, then for all simple objects $S \in \mathcal{A}$, there is at most one $T$ with $\Ext^1(S, T) \neq 0$.
	\end{enumerate}
\end{l1}

\begin{proof}
	We begin with the proof of (a). If (C1) is not fulfilled, then there is a finite set $\mathcal{T}$ of pairwise non-isomorphic simple objects with $S \in \mathcal{T}$ and
	\begin{equation*}
	d := \sum_{T \in \mathcal{T}} d^1_T(S, T) \ge 2.
	\end{equation*}
	We can assume that $d < \infty$ and the subcategory $\mathcal{A}'$ of objects of Loewy length $\le 2$ with composition factors in $\mathcal{T}$ is equivalent to a module category.
	
	In $\mathcal{A}'$, there is an indecomposable projective envelope $P_2(S)$ of $S$ with a socle of length $d$.
	
	Thus $P_2(S) \in \mathcal{A}' \subset \mathcal{A}$ is not colocal and by definition, $\mathcal{A}$ is not of colocal type.\\
	
	To prove (b), we assume that we can choose $\mathcal{T} = \{T_1, T_2, S\}$ so that $T_1 \ncong T_2$ and $\Ext^1(S, T_i) \neq 0$ for $1 \le i \le 2$. Then $T_1 \oplus T_2 \mid \soc P_2(S)$ and by Lemma \ref{LemSoc}, the lattice $\mathsf{S}(\mathcal{A})$ is not distributive.
\end{proof}

To show that (C2) holds if $\mathcal{A}$ is of colocal type, we use the following auxiliary lemma:

\begin{l1} \label{LemCommDiag}
	Let $\mathcal{A}$ be an abelian category and $S, T_1, \dots, T_n$ simple objects in $\mathcal{A}$ so that there are exact sequences with indecomposable middle terms
		\begin{equation*}
	\eta_i: \xymatrix{
		0 \ar[r] & S \ar^{f_i}[r] & X_i \ar[r] & T_i \ar[r]& 0
	}
	\end{equation*}
	for $1 \le i \le n$. Furthermore, suppose that $2 \le n$ and for all $1 \le i, j \le n$ either $T_i = T_j$ or $T_i \ncong T_j$ holds. For $f := \left[\begin{smallmatrix}
	f_1 \\f_2 \\ \vdots \\ f_{n-1} 
	\end{smallmatrix}\right]$ let there be a commutative diagram
			\begin{equation} \label{DiagComm}
	\xymatrix{ 
		S \ar^-{f}[r] \ar_{f_n}[d] & \bigoplus_{i = 1}^{n-1} X_i \ar^{g}[d]\\
		X_{n} \ar_-{g_n}[r] & X_n 
	}
	\end{equation}
	so that $g$ is an epimorphism and $g_n$ is an isomorphism.

	Then there are $1 \le i_1, \dots, i_m \le n-1$ so that $T_{i_1} = T_n, \dots, T_{i_m} = T_n$ and $\eta_n$ is linearly dependent of $\eta_{i_1}, \dots, \eta_{i_{m}}$ over $T_n$.
\end{l1}

\begin{proof}
	  Let $i: S^{n-2} \hookrightarrow \Ker g$ be the natural injection and $F := diag(f_1, f_2, \dots, f_{n-1})$.
	  
	   Then the following diagram is exact, because all rows and the first and second column are exact:
\begin{equation*}
\xymatrix@M=0.5pc@R=1.7pc@C=1.7pc{	&0\ar[d]&0\ar[d]&0\ar[d]&\\
	0 \ar[r] & S \ar_{\left[\begin{smallmatrix}
		\text{id} &\dots &\text{id}
		\end{smallmatrix}\right]}[d] \ar^{f_n}[r] &  X_n\ar_{gg_n^{-1}}[d] \ar[r] &  T_n \ar[r] \ar[d] & 0\\
	0 \ar[r] & S^{n-1} \ar^-{F}[r] \ar[d]&\bigoplus_{i = 1}^{n-1} X_i \ar[r]\ar[d] & \bigoplus_{i = 1}^{n-1} T_i \ar[d]\ar[r]& 0 \\
	0 \ar[r] & S^{n-2} \ar[d]  \ar@{^{(}->}^-{i}[r]& \Coker gg_n^{-1} \ar[d] \ar[r] & \Coker i \ar[r] \ar[d] & 0\\
	&0&0&0&}.
\end{equation*}
The second row of the diagram is $\bigoplus_{i=1}^{n-1}\eta_i$, while the first row of this diagram is $\eta_n$.

Thus, there are $1 \le i_1, \dots, i_m \le n-1$ so that $T_{i_1} = T_n, \dots, T_{i_m} = T_n$ and $\eta_n$ is linearly dependent of $\eta_{i_1}, \dots, \eta_{i_{m}}$ over $T_n$.
\end{proof}

Now we can show the following:

\begin{l1} \label{LemFourSimples1}
	Let $\mathcal{A}$ be an abelian length category.	
	\begin{enumerate}[(a)]
		\item If $\mathcal{A}$ is of colocal type, then (C2) is fulfilled.
		
		\item If $\mathsf{S}(\mathcal{A})$ is distributive, then for all simple objects $S \in \mathcal{A}$, there are no more than two non-isomorphic simple objects $T$ with $\Ext^1 (T, S) \neq 0$.
	\end{enumerate}
\end{l1}

\begin{proof}
	We start with the proof of (a). If (C2) is not fulfilled, then there is a simple object $S$ so that
	\begin{equation*}
	\sum_{T \text { simple}} d^1_T(T, S) \ge 3.
	\end{equation*}
	So there are three exact sequences
	\begin{equation*}
	\eta_i: \xymatrix{
		0 \ar[r] & S \ar^{f_i}[r] & X_i \ar[r] & T_i \ar[r]& 0
	}
	\end{equation*}
	for some indecomposable objects $X_i \in \mathcal{A}$ and simple $T_i$ with $i \in {1,2,3}$. If $T_i \cong T_j$ for some $i \neq j$, then we can assume that $T_i = T_j$. Furthermore, over $\End(T_i)^{op}$, $\eta_i$ is not a linear combination of the other two exact sequences.
	
	We consider the exact sequence
	\begin{equation*}
	\xymatrix@C=3pc{0 \ar[r] & S \ar^-{\left[\begin{smallmatrix}
			f_1 \\f_2\\f_3
			\end{smallmatrix}\right]}[r] & X_1 \oplus X_2 \oplus X_3 \ar^-{g}[r] & Y \ar[r] & 0}.
	\end{equation*}	
Since the objects $X_i$ are colocal, there is some indecomposable direct summand $Y'$ of $Y$ so that the morphism $g_1: X_1 \rightarrow Y'$ induced by $g$ is a monomorphism. 

If $g_i: X_i \rightarrow Y'$ is also induced by $g$ for $2 \le i \le 3$, then $\sum_{i=1}^3g_if_i = 0$. Thus there is at least one $i \in \{2,3\}$ so that $g_i$ is a monomorphism. By Lemma \ref{LemCommDiag}, this means $l(Y') \ge 3$. 

The same argument holds for every other direct summand $Y''$ of $Y$ so that the morphism $X_i \rightarrow Y''$ induced by $g$ is a monomorphism for some $1 \le i \le 3$. 

Since $l(Y) = 5$, $Y'$ is the only direct summand of $Y$ for which such a monomorphism exists. On the other hand, for $1 \le i \neq j \le 3$, $g$ induces a monomorphism $g': X_i \oplus X_j \rightarrowtail Y$.  With the arguments above, the image of the concatenation $S^2 \hookrightarrow X_i \oplus X_j \rightarrowtail Y$ is a subobject of $Y'$, so the kernel of $g'$ must be zero. Thus, $\mathcal{A}$ is not of colocal type.\\
	
	To prove (b), we assume that $T_1, T_2$ and $T_3$ are pairwise non-isomorphic.	
	
	Thus, we can assume without loss of generality that $T_1 \ncong S \ncong T_2$. Furthermore, 
	\begin{equation*}
	\begin{bmatrix}
	g'_1 & g'_2
	\end{bmatrix}:X_1 \oplus X_2 \rightarrowtail Y'
	\end{equation*}
	is a monomorphism. So both $T_1$ and $T_2$ arise only once as composition factors of $Y'$.
	
	 By Lemma \ref{LemTwoMonos}, there is a monomorphism $Y' \rightarrowtail \Coker g'_1 \oplus \Coker g'_2$. Thus $T_1$ is a composition factor of $Y'$, but not of $\Coker g'_1$ and $T_2$ is a composition factor of $Y'$, but not of $\Coker g'_2$. So $Y' \notin \sub \Coker g'_1$ and $Y' \notin \sub \Coker g'_2$. By Proposition \ref{PropDistributive}, $\mathsf{S}(\mathcal{A})$ is not distributive.
\end{proof}

\section{The third condition} \label{SecThird}

In this section, we prove that (C3) holds for all simple $S \in \mathcal{A}$ if $\mathcal{A}$ is of colocal type:

	\begin{enumerate}[\textup{(C3)}]
	\item If there is a simple object $S'$ with $\Ext^1(S,S') \neq 0$, let 
		\begin{equation*}
	\mathcal{T} := \{T \text{ simple and }\Ext^1(T, S) \neq 0  \mid \exists Z: l(Z) = 3, \soc Z = S',  \top Z = T\}
	\end{equation*}
	Then
	\begin{equation*}
	\sum_{T \in \mathcal{T}} d^1_T(T, S) \le 1.
	\end{equation*}
\end{enumerate}

Furthermore, we give some equivalent ways to define $\mathcal{T}$.

Proving that (C3) holds if $\mathcal{A}$ is of colocal type is more difficult than proving the other conditions. We need some auxiliary lemmas first:

\begin{l1} \label{Lem1}
	Let $Z_1, Z_2$ be objects of length 3 with socle $S'$ which are both local and colocal. Then for any monomorphisms $f_i: S' \rightarrowtail Z_i$, either the object $ Y := \Coker \left[\begin{smallmatrix}
	f_1 \\f_2
	\end{smallmatrix}\right]$ is indecomposable or there an isomorphism $\phi: Z_1 \rightarrow Z_2$ with $f_2 = \phi f_1$.
	\end{l1}

\begin{proof}
	Since $Z_i$ is colocal, $X_i := \Coker f_i$ is an indecomposable object of length 2.
	
	So there are epimorphisms $\left[\begin{smallmatrix}
	g_1 & g_2
	\end{smallmatrix}\right]: Z_1 \oplus Z_2 \twoheadrightarrow Y$ and $h_j: Z_j \twoheadrightarrow X_j$ so that the diagram below is exact for $1 \le i \neq j \le 2$, since its rows and the first two columns are exact:
	\begin{equation} \label{Diag}
	\xymatrix@M=0.5pc@R=1.7pc{&&0\ar[d]&0\ar[d]\\
		&0 \ar[r] \ar[d] &Z_i  \ar@{=}[r] \ar@{^{(}->}[d]& Z_i  \ar[r] \ar^{g_i}[d]& 0\\
		0 \ar[r] & S' \ar@{=}[d] \ar^-{\left[\begin{smallmatrix}
			f_1 \\f_2
			\end{smallmatrix}\right]}[r] & Z_1 \oplus Z_2 \ar[d] \ar^-{\left[\begin{smallmatrix}
			g_1 & g_2
			\end{smallmatrix}\right]}[r] & Y \ar[d] \ar[r] & 0\\
		0 \ar[r]	& S' \ar^{f_j}[r]\ar[d] & Z_j \ar^{h_j}[r]\ar[d] & X_j \ar[d]\ar[r] & 0 \\
		&0&0&0&}.
	\end{equation}
	Since $\soc Z_1 = S'$, there is some indecomposable direct summand $Y_1$ of $Y$ with a monomorphism $Z_1 \rightarrowtail Y_1$. We have
	\begin{equation*}
	5 = l(Y) \ge l(Y_1) \ge l(Z_1) = 3.
	\end{equation*}
	Thus $Y_1$ is the only direct summand of $Y$ with such a monomorphism. Analogously, there is a monomorphism $Z_2 \rightarrowtail Y_1$.
	
	Either $Y \cong Y_1$ and $Y$ is indecomposable or $Y = Y_1 \oplus Y_2$, where $1 \le l(Y_2) \le 2$. 
	
	Since $l(\top (Z_1 \oplus Z_2)) = 2$ and $\begin{bmatrix}
	g_1 & g_2
	\end{bmatrix}$ induces an epimorphism
	\begin{equation*}
	Z_1 \oplus Z_2 \twoheadrightarrow \top Y_1 \oplus \top Y_2,
	\end{equation*}
	both $Y_1$ and $Y_2$ are local. So for at least one $i \in \{1,2\}$, $g_i$ induces an epimorphism $Z_i \twoheadrightarrow \top Y_1$. 
	
	We write $g_i =: \left[\begin{smallmatrix}
	k_1 \\ k_2
	\end{smallmatrix}\right]: Z_i \rightarrowtail Y_1 \oplus Y_2$ and with (\ref{Diag}), there is an epimorphism $\left[\begin{smallmatrix}
	l_1 & l_2
	\end{smallmatrix}\right]: Y_1 \oplus Y_2 \twoheadrightarrow X_j$ with $j \neq i$ so that $l_1k_1 + l_2k_2 = 0$.

	If $l_1$ is an epimorphism, then $k_1l_1$ is non-zero; in particular it induces an isomorphism $\top Y_1 \rightarrow \top X_j$. So $\Coker k_1l_1 = 0$ and $k_1l_1$ is an epimorphism.
	
	It follows that $k_2l_2 = - k_1l_1$ is also an epimorphism. In particular, $l_2$ is an epimorphism. 		
	
	Since $X_j$ is local, at least one of $l_1$ and $l_2$ is an epimorphism  by Lemma \ref{LemTwoMonos}. So $l_2$ is always an epimorphism.
	
	Because of $l(Y_2) \le 2$, we see that $l_2$ is an isomorphism and
	\begin{equation*}
	l(Y_1) = l(Y) - l(Y_2) = 3.
	\end{equation*} 
	Since $g_1f_1 = - g_2f_2$ and $g_1$ and $g_2$ induce isomorphisms $Z_1 \cong Y_1 \cong Z_2$, there is an isomorphism $\phi: Z_1 \rightarrow Z_2$ so that $f_2 = f_1 \phi$.
\end{proof}

Let $\mathcal{T}$ be as in the definition of (C3) in the beginning of this section.

\begin{l1} \label{Lem2}
	Let $\mathcal{A}$ be an abelian length category. Suppose that there are simple objects $S, S'$ so that
	\begin{enumerate}
		\item $\Ext^1(S,S') \neq 0$
		\item $\sum_{T \in \mathcal{T}} d^1_T(T, S) \ge 2$
	\item $d^1_{S'}(S, S') = 1$ or $\mathcal{T}$ has two or more (non-isomorphic) elements.
	\end{enumerate}
	Then for some objects $Z_1$ and $Z_2$ of length 3 with socle $S'$ and top in $\mathcal{T}$, there are monomorphisms $f_i : S' \rightarrowtail Z_i$ for $1 \le i \le 2$ so that $\Coker \left[\begin{smallmatrix}
		f_1 \\f_2
	\end{smallmatrix}\right]$ is indecomposable.
\end{l1}

\begin{proof}
	Under these assumptions, there are $T_i \in \mathcal{T}$  for $1 \le i \le 2$ with indecomposable objects $X_i$ with exact sequences
	\begin{equation} \label{EqLinIndep}
	\eta_i: \xymatrix{0 \ar[r] & S \ar^{g_i}[r] & X_i\ar^{h_i}[r] & T_i \ar[r] & 0}
	\end{equation}
	that are not linearly dependent of each other over $\End(T_1)^{op}$. Furthermore, there are objects $Z_i$ for $1 \le i \le 2$ with exact sequences
		\begin{equation*} 
		\xymatrix{
			0 \ar[r] & S' \ar^{g'_i}[r] & Z_i\ar^{h'_i}[r] & X_i \ar[r] & 0}.
		\end{equation*}
	We begin by defining $f_1$ and $f_2$. Then we show that $Y := \Coker \left[\begin{smallmatrix}
		f_1 \\f_2
		\end{smallmatrix}\right]$ is indecomposable.
		
		If $Z_1 \ncong Z_2$, then we can choose arbitrary monomorphisms $f_i: S' \rightarrowtail Z_i$ for $1 \le i \le 2$.	Otherwise, we need to be more careful:
	
		If $Z_1 \cong Z_2$, we can assume that $Z_1 = Z_2$ and $T_1 = T_2$. Since $S' = \soc Z_1$, this implies $\Im g'_1 = \Im g'_2$. So $X_1 \cong X_2$.

		Let $X_3$ be the kernel of $Z_1 \twoheadrightarrow T_1$. Then there is an exact sequence
	\begin{equation*}
	\xymatrix{0 \ar[r] & S' \ar^{f}[r] & X_3 \ar^{g}[r] & S\ar[r] & 0}
	\end{equation*}
	and a monomorphism $f': X_3 \rightarrowtail Z_1$.
	
	We have $ \Im h'_if' \cong S \cong \Im g_ig$. By the assumptions, $d^1_{S'}(S, S') = 1$ and there is some isomorphism $\phi_i$ on $X_3$ so that $h'_if'\phi_i = g_ig$.
	
	Thus, we can define $f_i:= f'\phi_if$.
By \ref{Lem1}, $Y= \Coker \left[\begin{smallmatrix}
		f_1 \\f_2
		\end{smallmatrix}\right]$ is either indecomposable or there is some isomorphism $\phi: Z_1 \rightarrow Z_2$ with $f_2 = \phi f_1$.

	In the second case, there is an isomorphism $\psi: X_1 \rightarrow X_2$ and a commutative diagram with exact rows
		\begin{equation*}
		\xymatrix{0 \ar[r] & S' \ar^{f_1}[r] \ar@{=}[d]& Z_1 \ar^{\phi}[d] \ar^{h'_1}[r] & X_1 \ar^{\psi}[d]\ar[r] & 0\\
		0 \ar[r] & 	S' \ar^{f_2}[r] & Z_2  \ar^{h'_2}[r] & X_2\ar[r] & 0}
		\end{equation*}
	Since $f_i = f'\phi_if$, this induces monomorphisms $m_1$, $m_2$ with a commutative diagram
	\begin{equation*}
	\xymatrix{X_3 \ar^g[rr] \ar_{f'\phi_1}[dd] \ar^{f'\phi_2}[dr]& & X_3/S' = S \ar^<<<<<<{m_1}[dd] \ar^{m_2}[dr]\\
	& Z_2 \ar^<<<<<<<<<{h'_2}[rr] & & Z_2/S' = X_2\\ 
Z_1 \ar^{h'_1}[rr] \ar^{\phi}[ur] 	& & Z_1/S = X_1\ar_{\psi}[ur]}.
	\end{equation*}
	By definition, $h'_if'\phi_i = g_ig = m_ig$. Since $g$ is an epimorphism and $g_i$ a monomorphism, this implies $m_i = g_i$.
	
	So there is an epimorphism $\chi$ which makes the following diagram commutative with exact rows:
			\begin{equation*}
			\xymatrix{0 \ar[r] & S \ar^{g_1}[r] \ar@{=}[d]& X_1 \ar^{\psi}[d] \ar^{h_1}[r] & X_1 \ar^{\chi}[d]\ar[r] & 0\\
			0 \ar[r] & 	S \ar^{g_2}[r] & X_2  \ar^{h_2}[r] & X_2\ar[r] & 0}.
			\end{equation*}
		Thus, $\eta_1$ and $\eta_2$ are linearly dependent of each other over $\End(T_1)^{op}$, contrary to the assumptions.
\end{proof}

Now we can show the following:

\begin{l1} \label{LemFourSimples2}
	Let $\mathcal{A}$ be an abelian length category.	
	\begin{enumerate}[(a)]
		\item 	If $\mathcal{A}$ is of colocal type, then (C3) is fulfilled.
		
		\item Let $\mathsf{S}(\mathcal{A})$ be distributive with simple objects $S' \ncong S$ with $d_{S'}(S, S') = 1$. Then $\sum_{T \in \mathcal{T}} d_T(T, S) \le 1$.
	\end{enumerate}
\end{l1}

\begin{proof} 
	Let $S, S'$ be simple objects with  $d_{S'}(S, S') = 1$. By Lemma \ref{Lem2}, if we have $\sum_{T \in \mathcal{T}} d_T(T, S) \ge 2$, then there are some objects $Z_1$ and $Z_2$ of length 3 with socle $S'$ and tops $T_1, T_2 \in \mathcal{T}$, there are monomorphisms $f_i: S' \rightarrow Z_i$ so that $Y:= \Coker \left[\begin{smallmatrix}
	f_1 \\f_2
	\end{smallmatrix}\right]$ is indecomposable.

Since $d_{S'}(S, S') = 1$, there is up to isomorphism only one indecomposable object $X_3$ with an exact sequence of the form
	\begin{equation}
\xymatrix{0  \ar[r] &S' \ar^f[r] & X_3 \ar[r] & S \ar[r] &0}.
\end{equation}
So $X_3$ must be the kernel of $Z_1 \twoheadrightarrow T_1$ and $Z_2 \twoheadrightarrow T_2$. 

Furthermore, every monomorphism $S' \rightarrowtail X_3$ factors through $f$.  So $f_1$ and $f_2$ factor through $f$.
	
	Thus, the following diagram is commutative and exact, since its columns and the first two rows are exact:
	\begin{equation} \label{DiagCokerK}
	\xymatrix{&0\ar[d]&0 \ar[d]&\\		
		0 \ar[r] & S'  \ar@{=}[r] \ar^{\left[\begin{smallmatrix}
			f \\ f
			\end{smallmatrix}\right]}[d]& S'  \ar[r] \ar^{\left[\begin{smallmatrix}
			f_1 \\ f_2
			\end{smallmatrix}\right]}[d]& 0 \ar[d] \\
		0 \ar[r] & X_3 \oplus X_3 \ar[d] \ar[r] & Z_1 \oplus Z_2 \ar[d] \ar[r] & T_1 \oplus T_2 \ar@{=}[d] \ar[r] & 0\\
		0 \ar[r]	& S \oplus X_3  \ar[r]\ar[d] & Y  \ar[r]\ar[d] & T_1 \oplus T_2  \ar[d]\ar[r] & 0 \\
		&0&0&0&}.
	\end{equation}
	So $\soc  Y = S \oplus S'$.
	
	To show (a), we suppose that $\mathcal{A}$ is of colocal type, but (C3) is not fulfilled for some simple object $S \in \mathcal{A}$. 
	
	Then there is some $S'$ with $\Ext^1(S, S') \neq 0$ and we have $d_{S'}(S, S') = 1$ by Lemma \ref{LemTwoSimples}. With the arguments above, there is some indecomposable object $Y$ with  $\soc  Y = S \oplus S'$ and $\mathcal{A}$ is not of colocal type.\\
	
	Under the assumptions of (b), $d_{S'}(S, S') = 1$ and $S \ncong S'$. If $\sum_{T \in \mathcal{T}} d_T(T, S) \ge 2$, then by the arguments above, there is some indecomposable object $Y$ with $\soc  Y = S \oplus S'$. By Lemma \ref{LemSoc}, $\mathsf{S}(\mathcal{A})$ is not distributive.
\end{proof}

Under certain conditions, the definition of $\mathcal{T}$ is much simpler. To show this, we need some auxiliary lemmas. The first one will also be important in the next section, where we prove that every abelian length category that fulfils (C1) - (C3) is of colocal type.

\begin{l1} \label{LemLL2}
	Let $Y \in \mathcal{A}$ be an object with simple socle $S$, Loewy length 2 and a top of length $m$. 
	
	Then for any indecomposable object $X_m$ of length 2 with monomorphisms $f_m: S \rightarrowtail X_m$ and $f'_m: X_m \rightarrowtail Y$, there are indecomposable objects $X_1, \dots, X_{m-1}$ of length 2 with monomorphisms $f_i: S \rightarrowtail X_i$, $f'_i: X_i \rightarrowtail Y_i$ so that 
	\begin{equation*}
	f'_if_i = f'_mf_m
	\end{equation*}
	for all $1 \le i \le m$. They can be chosen so that over $\End^{op}(\Coker f_i)$ no exact sequence with monomorphism $f_i$ is a linear combination of exact sequences with monomorphisms in 
	\begin{equation*}
	\{f_1, \dots, f_{i-1}, f_{i+1}, \dots, f_m\}.
	\end{equation*}
\end{l1}

\begin{proof}
	We prove the lemma by induction on $m$: 
	
	It is obvious for $l(\top Y) = 1$ and we can assume that it has been proved for all objects $Y$ with Loewy length 2 and $l(\top Y) \le m$.
	
	So let $\top Y = \bigoplus_{i=1}^{m+1}T_i$ for some simple objects $T_i$. The kernel of $Y \twoheadrightarrow \bigoplus_{i=1}^{m}T_i$ is of length 2 and has the socle $S$; the kernel $Y'$ of $Y \twoheadrightarrow T_{m+1}$ has socle $S$, Loewy length 2 and $l(\top Y') = m$. 
	
	By Lemma \ref{LemTwoMonos}, there is an epimorphism
	\begin{equation*}
	\begin{bmatrix}
	f'_{m+1} &f'
	\end{bmatrix}: X_{m+1} \oplus Y' \twoheadrightarrow Y.
	\end{equation*}
	Since $l(X_{m+1} \oplus Y') = l(Y)+1$, the kernel of this epimorphism is $S$. So there are morphisms $f_{m+1}: S \rightarrow X_{m+1}$ and $f: S \rightarrow Y'$ with an exact sequence
	\begin{equation} \label{EqInduction}
	\xymatrix@C=4pc{0 \ar[r] & S \ar^-{\left[\begin{smallmatrix}
			f_{m+1} \\f
			\end{smallmatrix}\right]}[r] &X_{m+1} \oplus  Y' \ar^-{\left[\begin{smallmatrix}
			f'_{m+1} &f'
			\end{smallmatrix}\right]}[r] & Y\ar[r] &0 }.
	\end{equation}
	Since $l(Y') \ge 2$, there is some indecomposable object $X_m$ of length 2 with monomorphisms $f_m: S \rightarrowtail X_m$ and $f''_m: X_m \rightarrowtail Y'$ so that $f = f''_mf_m$. By the inductive assumption, we get objects $X_i$ of length 2 and monomorphisms $f_i: S \rightarrowtail X_i$ and $f''_i: X_i \rightarrowtail Y'$ so that $f''_if_i = f''_mf_m$ for $1 \le i \le m$. 
	
	With (\ref{EqInduction}), we have
	\begin{equation*}
	f'_{m+1}f_{m+1} = -f'f = -f'f''_mf_m =  -f'f''_if_i
	\end{equation*}
	for all $1 \le i \le m$.
	
	So we can set $f'_i := -f'f''_i$ and it remains to show that over $\End^{op}(\Coker f_i)$ there is no exact sequence with monomorphism $f_{m+1}$ that is a linear combination of exact sequences with monomorphisms in $\{f_1, \dots, f_m\}$:
	
	Otherwise, by definition of the Baer sum (see the beginning of Section \ref{SecExt}), there is some object $Z$ so that for $F = diag(f_1, \dots, f_m)$, there is the following commutative diagram with exact rows:
	\begin{equation*}
	\xymatrix{
		0 \ar[r] &	S \ar_-{f_{m+1}}[r] & X_{m+1} \ar[r] & T \ar[r] &0\\
		0 \ar[r] &		S^m \ar[r] \ar@{=}[d] \ar[u]& Z \ar[d] \ar[u] \ar[r] & T \ar[r] \ar@{=}[u] \ar[d] &0\\
		0 \ar[r] &		S^m \ar^-F[r] &\bigoplus_{i= 1}^m X_i\ar[r] & T^m \ar[r] &0}.
	\end{equation*}
	So $X_{m+1} \cong Z/S^{m-1}$ and the concatenation
		\begin{equation*}
	\xymatrix@C=4pc{ Z \ar[r] & \bigoplus_{i= 1}^m X_i\ar^-{\left[\begin{smallmatrix}
			f''_1 & \dots & f''_m
			\end{smallmatrix}\right]}[r] & Y'}
	\end{equation*}
	 induces a concatenation	
\begin{equation*}
f''_{m+1}: X_{m+1}\rightarrowtail \left(\bigoplus_{i= 1}^m X_i\right)/S^{m-1} \rightarrow Y'
\end{equation*}
	
	Since $f''_if_i = f''_mf_m$ for all $1 \le i \le m$, we get
	\begin{equation*}
	f''_{m+1} f_{m+1} = f''_mf_m = f.
	\end{equation*}
	But then $\Coker \left[\begin{smallmatrix}
	f_{m+1} \\f
	\end{smallmatrix}\right] \cong Y' \oplus T_{m+1}$, a contradiction to (\ref{EqInduction}).
\end{proof}

Analogously, the dual result holds. In particular, we get the following:

\begin{c1} \label{CorLL2}
	If there is an object $Y \in \mathcal{A}$ with simple socle $S$, Loewy length 2 and a top of length $m$, then $\sum_{T \text { simple}} d^1_T(T, S) \ge m$.
	
	Dually,  if there is an object $Y \in \mathcal{A}$ with simple top $S$, Loewy length 2 and a socle of length $m$, then $\sum_{T \text { simple}} d^1_T(S, T) \ge m$.
\end{c1}

As with modules, we call an object \textit{uniserial} if it has a unique composition series.

With the corollary above, we can show the following Lemma, which we will use to show the first two simplifications of $\mathcal{T}$:

\begin{l1} \label{LemUniserial2}
	If (C1) holds for an abelian length category $\mathcal{A}$, then every object in $\mathcal{A}$ that is both local and colocal is uniserial.
\end{l1}

\begin{proof}	
	Suppose that $X$ is an object that is both local and colocal. Then every non-zero factor object of $X$ is indecomposable. Set $X_0 := X$ and choose $X_1, \dots, X_n$ with epimorphisms
	\begin{equation*}
	X_0 \twoheadrightarrow X_1 \twoheadrightarrow \dots \twoheadrightarrow X_n = 0
	\end{equation*}
	so that $\Ker(X_{i-1} \twoheadrightarrow X_i)$ is simple for all $1 \le i \le n$.
	
	If $X$ is not uniserial, then it has a factor object which is not colocal. In particular, there is some maximal integer $m$, so that $X_m$ is not colocal. For every $S \mid \soc X_m$, the quotient $X_m / S$ is uniserial and thus $l(\soc X_m) = 2$. 
	
	The kernel $X'$ of $X_m \twoheadrightarrow X_{m+3}$ has the length 3 and $\soc X' = \soc X_m$. So the Loewy length of $X'$ is 2 and by Corollary \ref{CorLL2} and (C1), $l(\top X') \ge 2$. Thus $l(\soc X') + l(\top X') > l(X')$ and there is some simple object $S \mid X'$. 
	
	For $S'$ with $\soc X' = S \oplus S'$, the quotient $X'/S'$ is decomposable. But $X'/S'$ is the kernel of the epimorphism $X_m/S' \twoheadrightarrow X_{m+3}$, which contradicts the assumption that $X_m/S$ is uniserial.
	
	So if  $\mathcal{A}$ fulfils (C1), then every object $X$ that is both local and colocal is also uniserial.
\end{proof}

We still need two very different auxiliary lemmas, which we will use to show the final simplification of $\mathcal{T}$:

Recall that the $\Ext$-quiver of $\mathcal{A}$ is a quiver that has the simple objects of $\mathcal{A}$ as vertices and an arrow between vertices $S$ and $T$ if and only if $\Ext^1(S, T) \neq 0$.

\begin{l1} \label{LemAux2}
	Let $\mathcal{A}$ be an abelian length category that fulfils (C1).
	
	Suppose that there are simple objects $S', T$ so that no cycle in the $\Ext$-quiver has $S'$ as a vertex and there is an indecomposable object $Z$ of length 3 with $\soc Z = S'$ and $\top Z = T$.

	Then $\Ext^2(T, S') = 0$. 

\end{l1}

\begin{proof}
	By the assumptions on $Z$, there is a simple object $S$ with
	\begin{equation*}
	\Ext^1(S, S') \neq 0 \neq \Ext^1(T, S).
	\end{equation*}
	Since no cycle in the $\Ext$-quiver has $S'$ as a vertex,
	\begin{equation} \label{EqExt=0}
	\Ext^1(S', S')  = \Ext^2(S',S')= 0.
	\end{equation}
	In particular, $S' \ncong S$. Furthermore, an exact sequence
	\begin{equation*}
	\xymatrix{0 \ar[r] & S' \ar^f[r] & X \ar^g[r] & Y \ar^h[r] & Z \ar[r] & 0} \in \Ext^2(Z, S')
	\end{equation*}
	induces a commutative diagram with exact rows 
		\begin{equation*}
	\xymatrix{0 \ar[r] & S' \ar^f[r] \ar@{=}[d] & X \ar^g[r] \ar@{=}[d]& Y \ar^h[r]\ar^{m_1}[d] & Z \ar[r] \ar^{m_2}[d]& 0\\
		0 \ar[r] & S' \ar^f[r] & X \ar^{g'}[r] & Y' \ar^{h'}[r] & S' \ar[r] & 0}
	\end{equation*}
	with monomorphisms $m_1$ and $m_2$. Since $\Ext^2(S',S') = 0$, the second row of the diagram splits. Thus, the first row also splits and 
		\begin{equation} \label{EqExt=02}
	\Ext^2(Z, S') = 0.
	\end{equation}
	Since (C1) holds,
	\begin{equation} \label{EqC1}
	\sum_{T' \text{ simple}} d^1_{T'}(S, T') = 1 = \sum_{T' \text{ simple}} d^1_{T'}(T, T').
	\end{equation}
	So $\Ext^1(S, S) = 0$ and thus $S \ncong T$ and $\Ext^1(T, S') = 0$.
	
	Le $X$ be an indecomposable object with exact sequences	
	\begin{equation} \label{SeqZ}
	\xymatrix{0 \ar[r] & S' \ar[r] & Z \ar[r] & X \ar[r] & 0}
	\end{equation}
	and
		\begin{equation} \label{SeqX}
	\xymatrix{0 \ar[r] & S \ar[r] &  X \ar[r] & T\ar[r] & 0}.
	\end{equation}

	The short exact sequence (\ref{SeqZ}) induces a long exact sequence
\begin{equation*}
\xymatrix{ \dots \ar[r] &\Ext^1(S', S') \ar[r]	& \Ext^2(X, S') \ar[r] &\Ext^2(Z, S') \ar[r] & \dots}.
\end{equation*}
 So $\Ext^2(X, S') = 0$ by (\ref{EqExt=0}) and (\ref{EqExt=02}).
 
Furthermore, (\ref{SeqX}) induces a long exact sequence
 \begin{equation*}
\xymatrix{ \dots \ar[r] &\Ext^1(T, S') \ar[r] &\Ext^1(X, S')\ar^{\alpha}[r] &\Ext^1(S, S') 	\ar `r[r] `d[l] `l[dlll] `d[r] [dll]& \\
	 &\Ext^2(T, S') \ar[r] &\Ext^2(X, S')\ar[r] &\dots}.
 \end{equation*}

Since $\Ext^1(T, S') = 0$, the morphism $\alpha$ is a monomorphism. With (\ref{EqC1}),
\begin{equation*}
\dim_{\End(S')} \Ext^1(X, S')  \ge 1 = \dim_{\End(S')} \Ext(S, S'),
\end{equation*}
 so $\alpha$ must be an isomorphism. We get
\begin{equation*}
\Ext^2(T, S') \cong \Ext^2(X, S') = 0
\end{equation*}
and the proof is complete.
\end{proof}

The other direction of \ref{LemAux2} is even more generally true:

\begin{l1} \label{LemAux}
	Let $S, S', T \in \mathcal{A}$ be simple objects with $\Ext^1(S, S') \neq 0 \neq \Ext^1(T, S)$ so that $\Ext^1(T, S') = 0 = \Ext^2(T, S')$. Then there is an object $Z$ of length 3 so that $\soc Z = S'$ and $\top T = T$.  
\end{l1}

\begin{proof}
	By the assumptions, there is an indecomposable object $X$ with an exact sequence	
	\begin{equation*}
	\xymatrix{
		0 \ar[r] & S' \ar[r] & X \ar[r] & S \ar[r] & 0}.
	\end{equation*}
	This induces a long exact sequence
	\begin{equation*}
	\xymatrix{\dots \ar[r] & \Ext^1(T, S') \ar[r] &  \Ext^1(T, X) \ar[r] &  \Ext^1(T, S) \ar[r] & \Ext^2(T, S') \ar[r] &  \dots}.
	\end{equation*}
	Since $\Ext^1(T, S') = \Ext^2(T, S') = 0$, we see that $\Ext^1(T, X) \cong \Ext^1(T, S) \neq 0$. So there is some indecomposable object $Z$ with an exact sequence
		\begin{equation*}
	\xymatrix{
		0 \ar[r] & X \ar[r] &Z \ar[r] & T \ar[r] & 0}
	\end{equation*}
	and obviously $l(Z) =3$, $\soc Z = S'$ and $\top T = T$.
\end{proof}

The following equivalence holds:

\begin{p1} \label{PropEquiv}
	Suppose that $\mathcal{A}$ is an abelian length category which fulfils (C1).
	
	For fixed simple objects $S$ and $S'$ with $\Ext^1(S, S') \neq 0$, the class
	\begin{equation*}
	\mathcal{T} = \{T \text{ simple and }\Ext^1(T, S) \neq 0  \mid \exists Z: l(Z) = 3, \soc Z = S',  \top Z = T\}
	\end{equation*}
	is the same as
	\begin{equation*}
\mathcal{T}' := \{T \text{ simple and }\Ext^1(T, S) \neq 0 \mid \exists Z: l(Z) \ge 3, \soc Z = S',  \top Z = T\}.
\end{equation*}	
If $S \ncong S'$, then $\mathcal{T} = \mathcal{T}''$ with 
	\begin{equation*}
\mathcal{T}'' := \{T \text{ simple and }\Ext^1(T, S) \neq 0 \mid \exists Z: \soc Z = S',  \top Z = T\}.
\end{equation*}
If no cycle in the $\Ext$-quiver has $S'$ as a vertex, then $\mathcal{T} = \mathcal{T}'''$ with
\begin{equation*}
\mathcal{T}''' := \{T \text{ simple} \mid \Ext^1(T, S) \neq 0 = \Ext^2(T, S')\}.
\end{equation*}
\end{p1}

\begin{proof}
	The relations $\mathcal{T} \subset \mathcal{T}' \subset \mathcal{T}''$ are clear.
	
	 On the other hand, suppose that $\Ext^1(T, S) \neq 0$ and there is some object $Z$ with $l(Z) = 3$, $\soc Z = S'$ and $\top Z = T$.  We can define $Z_0 := Z$ and $Z_i := Z_{i-1}/\soc Z_{i-1}$ for $i \in \mathbb{N}$. Then there is some $Z_m$ with length $l(Z_m) = 3$. By Lemma \ref{LemUniserial2}, $Z$ and thus $Z_m$ is uniserial. Furthermore, $\top Z_m = T$ by construction.
	 
	 By (C1), $\Ext^1(T, T') = 0$ for all $T' \ncong S$ and $\Ext^1(S, T') \neq 0$ for all $T' \ncong S'$. Thus $\soc Z_m = S'$ and $\mathcal{T}' \subset \mathcal{T}$.
	
	Now suppose that $S \ncong S'$ and there is some simple object $T$ with $\Ext^1(T, S) \neq 0$. Since $\mathcal{A}$ fulfils (C1), we have $\Ext^1(T, S') = 0$. So every object $Z$ with $\soc Z = S'$ and $\top Z = T$ has length greater or equal to 3 and thus $\mathcal{T}' = \mathcal{T}''$.
	
	Finally, if no cycle in the $\Ext$-quiver has $S'$ as a vertex, then $\mathcal{T} = \mathcal{T}'''$ by Lemma \ref{LemAux2} and Lemma \ref{LemAux}.                                                                                 
\end{proof}

\section{An Equivalence Theorem} \label{SecEquiv}

In the last sections, we have shown that (C1) - (C3) has to be fulfilled if $\mathcal{A}$ is of colocal type. In this section, we prove the other direction.

Reformulating \cite{ASS}, Chapter V, Theorem 2.6, an Artin algebra $A$ is right serial (that is, every right projective module over $A$ is uniserial) if and only if (C1) holds.

Note that in this case a uniserial object is uniquely determined up to isomorphism by its composition series (see \cite{ASS}, Chapter V, 2.7); in fact, it is even uniquely determined up to isomorphism by its top and its length: 

If $T \in \mathcal{A}$ is simple, then there is exactly one maximal path in the $\Ext$-quiver of $\mathcal{A}$ that starts in $T$:
\begin{equation*}
T = T_1 \rightarrow T_2 \rightarrow T_3 \rightarrow \dots.
\end{equation*}
 Every uniserial object with top $T$ and length $n$ has the socle $T_n$. If we denote this object by $U_{T, n}$, then we have epimorphisms
\begin{equation*}
\dots \twoheadrightarrow U_{T, n} \twoheadrightarrow U_{T, n-1} \twoheadrightarrow \dots \twoheadrightarrow U_{T, 2} \twoheadrightarrow U_{T, 1} \cong T \twoheadrightarrow U_{T, 0} := 0.
\end{equation*}

On the other hand, suppose that (C2) and (C3) hold. Let
\begin{equation*}
\dots \rightarrow S_3 \rightarrow S_2 \rightarrow S_1 = S
\end{equation*}
be a maximal path in the $\Ext$-quiver of $\mathcal{A}$ so that there is an object $Z$ with socle $S$ and top $S_n$.

Then there is at most one different maximal path
\begin{equation*}
\dots \rightarrow S'_3 \rightarrow S'_2 \rightarrow S'_1 =S
\end{equation*}
in the $\Ext$-quiver so that an object with socle $S$ and top $S'_n$ exists. If such a path exists, then $\sum_{T \text { simple}} d^1_T(T, S) = 2$ and $S_1 \ncong S'_1$.

If there is only one such path but $\sum_{T \text { simple}} d^1_T(T, S) = 2$, then we set $S'_1 := S_1, \dots, S'_n := S_n$.

With this, we can formulate the next lemma:

\begin{l1} \label{LemTopL2}
		Suppose that $\mathcal{A}$ fulfils (C1) - (C3). Let $Y \in \mathcal{A}$ be a colocal object with $\soc Y =: S$ that is not local. Then $\sum_{T \text { simple}} d^1_T(T, S) = 2$, $l(\top Y) = 2$ and there are $m, m' \in \mathbb{N}$ so that $Y/S \cong U_{S_m, m-1} \oplus U_{S'_{m'}, m'-1}$.
\end{l1}

\begin{proof}
	Let $Y$ be a colocal object with socle $S$ and $\soc (Y/S) := \bigoplus_{i=1}^n R_i$,  with simple $S_i$ for $1 \le i \le n$ and $n \neq 0$. Then there is a subobject $Y'$ of $Y$ with $\soc Y' = S$ and $Y'/S =  \bigoplus_{i=1}^n R_i$. So $Y'$ has Loewy length 2 and by Corollary \ref{CorLL2}, (C2) implies that $n \le 2$ and $\{R_1, R_2\} = \{S_1, S'_1\}$ if $n = 2$.
	
	If $S_1 \mid \soc (Y/S)$, then there are simple objects $R'_k$, $1 \le k \le n'$ so that
	\begin{equation*}
	\soc (Y/S)/S_1 = S'_1 \oplus \bigoplus_{i = 1}^{n'} R'_k.
	\end{equation*}
	Thus, for $1 \le k \le n'$, there are subobjects $Z_k$ of $Y$ so that $\soc Z_k = S$, $\top Z_k = R'_k$ and $l(Z_k) \ge 3$. By Proposition \ref{PropEquiv}, (C3) implies that $n' \le 1$ and $R'_1 = S_2$ if it exists. 
	
	Inductively, there is some $U_{S_m, m}$ with a monomorphism $U_{S_m, m} \rightarrowtail Y$.
	
	The analogous argument holds for $S'_1$ and there is some $U_{S'_{m'}, m'}$ with a monomorphism $U_{S'_{m'}, m'}$.
	
	We can either choose these objects so that $\top Y = S_m \oplus S'_{m'}$ or $l(\top Y) = 1$.
	
	But the latter is not possible, since otherwise there is some subobject $Y''$ of a factor object of $Y$ which has a simple top $T$ and there are $m, m' \in \mathbb{N}$ so that $\Ker(Y'' \twoheadrightarrow T) = S'_{m'} \oplus S_{m}$. By Lemma \ref{LemLL2}, this contradicts (C1).
	
	So $\top Y = S_m \oplus S_{m'}$ and the monomorphisms $U_{S_{m}, m} \rightarrowtail Y$, $U_{S'_{m'}, m'} \rightarrowtail Y$ imply that $Y/S = U_{S_{m}, m-1} \oplus U_{S'_{m'}, m'-1}$.
\end{proof}

As a corollary, we get:

\begin{c1} \label{CorTopL2}
	Suppose that $\mathcal{A}$ fulfils (C1) - (C3). Let $Y \in \mathcal{A}$ be a colocal object with $\soc Y = S$ that is not local. Then there are $m, m' \ge 2$ with exact sequences
	\begin{equation} \label{Seq1}
	\xymatrix{0 \ar[r] & U_{S_m, m} \ar[r] & Y \ar[r] & U_{S'_{m'}, m'-1}}
	\end{equation}
	and
	\begin{equation} \label{Seq2}
	\xymatrix{0 \ar[r] & U_{S'_{m'}, m'} \ar[r] & Y \ar[r] & U_{S_m, m-1}}.
	\end{equation}
\end{c1}

Furthermore, we see:

\begin{l1} \label{LemLocal}
 If $\mathcal{A}$ fulfils (C1), then every object that is local is also colocal.
\end{l1}

\begin{proof}
	Let $Y$ be local but not colocal. Then every quotient of $Y$ is local. Inductively, we can assume that every real quotient of $Y$ is local and colocal and thus uniserial by Lemma \ref{LemUniserial2}. With $Y':= Y/soc Y$ and $Y'' := Y'/\soc Y'$, the object 
	\begin{equation*}
	X := \ker(Y \twoheadrightarrow Y'')
	\end{equation*} 
	has Loewy length 2. Suppose that
	\begin{equation*}
	T_1 \oplus T_2 \mid \soc Y = \soc X.
	\end{equation*}
	Then $X/T_i$ is a subobject of $Y/T_i$ for $1 \le i \le 2$ and thus uniserial. So $X$ is local and by Corollary \ref{CorLL2}, (C1) does not hold, a contradiction to the assumption. 
\end{proof}

With this, we can prove the following:

\begin{l1} \label{LemFactor}
	Let $\mathcal{A}$ be an abelian length category for which (C1) holds.
	
	If $g_m: U_{T, m} \twoheadrightarrow T$ and $g_n: U_{T, n} \twoheadrightarrow T$ are epimorphisms and $m \le n$, then $g_n$ factors through $g_m$.
\end{l1}

\begin{proof}
	In  this case, there is some morphism $g'_n:  U_{T, n} \twoheadrightarrow T$ which factors through $g_m$. We show by induction that there is an isomorphism $\phi$ on $U_{T, n}$ with $g_n = g'_n \phi$: 
	
	For $n=1$, this is clear and for $n = 2$, it follows from $\dim_{\End(T_2)}\Ext^1(T, T_2) = 1$.
	
	Suppose that the assertion is true for all morphisms $U_{T, n'} \twoheadrightarrow T$ with $1 \le n' < n$. 
	
	There is some object $Y$ with an exact sequence
	\begin{equation*}
	\xymatrix@C=3pc{0 \ar[r] & Y \ar[r] & U_{T, n} \oplus U_{T, n }\ar^-{\left[\begin{smallmatrix}
			g_n & g'_n
			\end{smallmatrix}\right]}[r] & U_{T, n-1} \ar[r] &0}.
	\end{equation*}
	Since $l(Y) = l( U_{T, n})+1$ and $T_n \mid \soc Y$, the object $Y$ is either local or of the form $U_{T, n} \oplus T_n$. 
	
	In the latter case, there is an isomorphism $\phi$ on $U_{T, n}$ with $g_n = g'_n \phi$.
	
	The former case is impossible by Corollary \ref{LemLocal}: Since $Y$ is not a subobject of $U_{T, n}$, it is not colocal by Lemma \ref{LemTwoMonos}.
\end{proof}

Before we can proof the main theorem, we still need one lemma about the objects in $\mathcal{A}$:

\begin{l1} \label{LemFactor2}
	Suppose that $\mathcal{A}$ fulfils (C1) - (C3). Let $Y, Y'$ be colocal objects with socle $S$ so that there is a simple object $T$ with $T \mid \top Y$ and $T \mid \top Y'$. If $l(Y) \le l(Y')$, then every epimorphism $g: Y \twoheadrightarrow T$ factors through every epimorphism $g': Y' \twoheadrightarrow T$.
\end{l1}

\begin{proof}
	It suffices to prove the case $l(Y) = l(Y')$: If $l(Y) < l(Y')$, then $Y'$ has a subobject $Y''$ with $l(Y'') = l(Y)$, $\soc Y'' = S$ and $T \mid \top Y''$. For a monomorphism $f: Y'' \rightarrowtail Y'$, we can define $g'' = g'f$ and if $g$ factors through $g''$, then it also factors through $g'$.
	
	If $Y$ is local, then this lemma is already the result of Lemma \ref{LemFactor}.
	
	To prove the assertions for non-local objects with $l(Y) = l(Y')$ we use Corollary \ref{CorTopL2}:	There are $m, m' \ge 2$ with an exact sequence
	\begin{equation*}
	\xymatrix{0 \ar[r] & U_{S_m, m} \ar[r] & Y \ar[r] & U_{S'_{m'}, m'-1}}.
	\end{equation*}
	Without loss of generality, we can assume that $T = S'_{m'}$ and it suffices to show that
	\begin{equation} \label{EqDim1}
	\dim_{\End(U_{S'_{m'}, m'-1})} \Ext^1(U_{S'_{m'}, m'-1}, U_{S_m, m}) = 1.
	\end{equation}
	The exact sequences
		\begin{equation*}
	\xymatrix{0 \ar[r] & S_i \ar[r] & U_{S_{m}, m-i+1} \ar[r] & U_{S_{m}, m-i} \ar[r] & 0}
	\end{equation*}
	for $1 \le i < m$
	induce long exact sequences
	\begin{equation*}
	\xymatrix{\dots \ar[r] &\Ext^1(U_{S'_{m'}, m'-1}, S_i)\ar[r] &\Ext^1(U_{S'_{m'}, m'-1}, U_{S_m, m-i+1})	\ar `r[r] `d[l] `l[dll] `d[r] [dl]& \\
	&	\Ext^1(U_{S'_{m'}, m'-1}, U_{S_m, m-i})\ar[r] &\dots ~.}
	\end{equation*}
	If $S \ncong S_i$ for $2 \le i \le m$, then $\Ext^1(S'_2, S_i) = 0$ and thus
	\begin{equation*}
	\Ext^1(U_{S'_{m'}, m'-1}, S_i) = 0.
	\end{equation*} 
	Since $U_{S_m, m-(m-1)} = S_i$, we get inductively that
	\begin{equation*}
	\Ext^1(U_{S'_{m'}, m'-1}, U_{S_m, m-i+1}) = 0
	\end{equation*} 
	for $2 \le i < m$. 
	
	Thus,
	\begin{equation*}
	\Ext^1(U_{S'_{m'}, m'-1}, U_{S_m, m}) \cong \Ext^1(U_{S'_{m'}, m'-1}, S_1).
	\end{equation*} 
	By Lemma \ref{LemFactor}, $\dim_{\End(U_{S'_{m'}, m'-1})} \Ext^1(U_{S'_{m'}, m'-1}, S_1)= 1$ and thus (\ref{EqDim1}) holds.
	
	Now suppose that there is some $2 \le i \le m$ with $S_i \cong S$. Then $d^1_{S_{j-1}}(S_j, S_{j-1}) =1$ for all $1 \le j \le i$. Since $S_1 \cong S \cong S_i$, we get $\End(S_j) \cong \End(S)$ for all $1 \le j \le i$.
	
	By Lemma \ref{LemTopL2}, $\sum_{T \text { simple}} d^1_T(T, S) = 2$ and by definition of $S_2, S'_2$ we get $S'_2 \ncong S_2$. With (C1), this implies $S'_k \ncong S_j$ for all $2 \le k \le m'$ and $1 \le j \le i$.
	
	In particular, $S'_k \ncong S$ for $2 \le k \le m'$.

	The exact sequence
\begin{equation*}
\xymatrix{0 \ar[r] & U_{S'_{m'-1}, m'-2} \ar[r] & U_{S'_{m'}, m'-1} \ar[r] & S_{m'} \ar[r] & 0}
\end{equation*}
induces a long exact sequence
\begin{equation} \label{ExLong}
\xymatrix{\dots \ar[r] &\Ext^1(S'_{m'}, U_{S_m, m})\ar[r] &\Ext^1(U_{S'_{m'}, m'-1}, U_{S_m, m})	\ar `r[r] `d[l] `l[dll] `d[r] [dl]& \\
	&	\Ext^1(U_{S'_{m'-1}, m'-2}, U_{S_m, m})\ar[r] &\dots}
\end{equation}	
If $l(U_{S'_{m'}m'-1}) > 1$, then $\Ext^1(S'_{m'}, U_{S_m, m})=0$ and there is a monomorphism 
\begin{equation*}
\Ext^1(U_{S'_{m'}, m'-1}, U_{S_m, m}) \rightarrowtail \Ext^1(U_{S'_{m'-1}, m'-2}, U_{S_m, m}).
\end{equation*}
So inductively, it suffices to show the lemma for $l(U_{S'_{m'},m'-1}) = 1$. In this case, $m' = 2$, $U_{S'_{m'-1}, m'-2} = 0$ and (\ref{ExLong}) yields an epimorphism
\begin{equation*}
\Ext^1(S'_{2}, U_{S_m, m}) \twoheadrightarrow \Ext^1(U_{S'_{2}, 1}, U_{S_m, m}).
\end{equation*}
The short exact sequence
\begin{equation*}
\xymatrix{0 \ar[r] &S \ar[r] & U_{S_{m}, m} \ar[r] & U_{S_{m}, m-1} \ar[r] &  0}
\end{equation*}
induces a long exact sequence
\begin{equation*}
\xymatrix{\dots \ar[r] &\Ext^1(S'_2, S) \ar[r] & \Ext^1(S'_{2}, U_{S_m, m}) \ar[r] &\Ext^1(S'_{2}, U_{S_m, m-1}) \ar[r] & \dots }.
\end{equation*}
Since $S'_2 \ncong S$, we have $\Ext^1(S'_{2}, U_{S_m, m-1}) = 0$ and an epimorphism
\begin{equation*}
\Ext^1(S'_2, S)\twoheadrightarrow \Ext^1(S'_{2}, U_{S_m, m}).
\end{equation*}
Since $S_2 \ncong S'_2$, we have $\End(S_i) \cong \End(S'_j) \cong \End(S)$ for all $1 \le i \le m$ and $1 \le i \le m'$. In particular, $\End(U_{S'_{m'}, m'-1}) \cong \End(S)$.

With the arguments above,
\begin{equation*}
1 = d^1_{S}(S'_2, S) = \dim_{\End(S)}\Ext^1(S'_2, S)
\end{equation*}
implies (\ref{EqDim1}).
\end{proof}

Now we can prove Theorem \ref{ThShape}:

\begin{t1} \label{ThShape2}
	An abelian length category $\mathcal{A}$ is of colocal type if and only if the following conditions hold for all simple objects $S \in \mathcal{A}$:

\begin{equation*} \tag{C1}
\sum_{T \text{ simple}} d^1_T(S, T) \le 1
\end{equation*}

\begin{equation*} \tag{C2}
\sum_{T \text { simple}} d^1_T(T, S) \le 2
\end{equation*}
\begin{enumerate}[\textup{(C3)}]
	\item If there is a simple object $S'$ with $\Ext^1(S,S') \neq 0$, let 	\begin{equation*}
	\mathcal{T} := \{T \text{ simple and }\Ext^1(T, S) \neq 0  \mid \exists Z: l(Z) = 3, \soc Z = S',  \top Z = T\}
	\end{equation*}
	Then
	\begin{equation*}
	\sum_{T \in \mathcal{T}} d^1_T(T, S) \le 1.
	\end{equation*}
\end{enumerate}
\end{t1}

\begin{proof}
Suppose that (C1) - (C3) holds, but $\mathcal{A}$ is not of colocal type. Then by Lemma \ref{LemUniserial}, there is an exact sequence
	\begin{equation} \label{EqSeq}
\xymatrix@C=4pc{0 \ar[r] & X \ar[r] & \oplus_{i=1}^m Y_i \ar^-{\left[\begin{smallmatrix}
	f_1 & \dots & f_m
	\end{smallmatrix}\right]}[r] & T \ar[r] & 0}.
\end{equation}
with $X$ indecomposable, $T$ simple, $Y_i$ colocal for $1 \le i \le m$ and $m \ge 2$. In particular, $f_i \neq 0$ and $T \mid \top Y_i$ for all $1 \le i \le m$.

We can order the objects $Y_i$ so that for some $m', n \in \mathbb{N}$, we have $\soc Y_i = T_n$ for $1 \le i \le m'$ and $\soc Y_i = T_{n_i}$ with $n_i > n$ for $m' < i \le m$. Furthermore, we can assume that $l(Y_1) \ge l(Y_2) \ge \dots \ge l(Y_{m'})$.

By Corollary \ref{CorTopL2} and Lemma \ref{LemFactor}, we can choose an epimorphism $g_n: U_{T, n} \twoheadrightarrow T$ so that $f_i$ factors through $g_n$ for $m' < i \le m$. On the other hand, by Lemma \ref{LemFactor2}, $f_i$ factors through $f_1$ for $1 \le i \le m'$.

It remains to show that $g_n$ factors through every epimorphism $f_1$. By Corollary \ref{CorTopL2}, there is a monomorphism $U_{T, n} \rightarrowtail Y_1$ and $f_1$ induces an epimorphism $g_{n-1}: U_{T, n-1} \twoheadrightarrow T$. By Lemma \ref{LemFactor}, $g_n$ factors through $g_{n-1}$. Thus $g_n$ also factors through $f_1$.

It follows that $f_i$ factors through $f_1$ for $1 \le i \le m$ and
	\begin{equation} \label{EqKer}
\Ker \begin{bmatrix}
f_1 & \dots & f_m
\end{bmatrix} \cong \bigoplus_{i=1}^{m-1} Y_i \oplus \Ker g_m
\end{equation}
and either $m = 1$ or $X$ is decomposable, contrary to the assumptions. 
\end{proof}

\section{The lattice $\mathsf{S}(\mathcal{A})$} \label{SecLattice}

We show in this section that the lattice $\mathsf{S}(\mathcal{A})$ is in fact the Cartesian product of certain sublattices. 

For Artin algebras $A$ of colocal type over algebraically closed fields, we will use this in Section \ref{SecStructure}, where we completely describe their lattice $\mathsf{S}(\mod A)$.

We begin with the following lemma:

\begin{l1} \label{LemSimpleObjects}
	Suppose $X$ is an indecomposable object and there is an index set $I$ with a monomorphism $X \rightarrowtail \bigoplus_{i \in I} Y_i.$ Set
	\begin{equation*}
	I' = \{ i \in I ~|~ \text{there is a simple object }S\text{ with } S \subseteq X \text{ and } S \subseteq X_i \}.
	\end{equation*}
	Then there is a monomorphism $X \rightarrowtail \bigoplus_{i \in I'} Y_i$.
\end{l1}

\begin{proof}
	There is a monomorphism
	\begin{equation*}
	\begin{bmatrix}
	f_1 \\ f_2
	\end{bmatrix}: X \rightarrowtail \bigoplus_{i \in I'} Y_i \oplus \bigoplus_{i \in I \setminus I'} Y_i
	\end{equation*}
	with
	\begin{equation*}
	f_1: X \rightarrow \bigoplus_{i \in I'} Y_i
 \quad \text{and} \quad
	f_2: X \rightarrow \bigoplus_{i \in I \setminus I'} Y_i.
	\end{equation*}
	With $i_1:\Ker (f_1) \rightarrowtail X$, the concatenation $f_2i_1$ must be a monomorphism. So there is no simple $S \subset \Ker (f_1)$ and thus $\Ker (f_1) = 0$, which implies that $f_1$ is a monomorphism.
\end{proof}

To simplify the notation, we define:

\begin{d1}
	For a class $\mathcal{M}$ of indecomposable objects in $\mathcal{A}$ let
	\begin{equation*}
\mathsf{S}(\mathcal{M}) := \mathsf{S}( \add \mathcal{M}).
	\end{equation*}
\end{d1}

\begin{l1} \label{LemSublattice}
	Let $\mathcal{M}$ be a class of indecomposable objects in $\mathcal{A}$. If
	\begin{equation} \label{EqSubinM}
	\ind \sub \mathcal{M} = \mathcal{M},
	\end{equation}
	then $\mathsf{S}(\mathcal{M})$ is a sublattice of $\mathsf{S}(\mathcal{A})$.
\end{l1}

\begin{proof}
	We need to show that for $C, C' \in \mathsf{S}(\mathcal{M})$, the join and the meet are again in $\mathsf{S}(\mathcal{M})$. 
	Since $C \wedge C' = C \cap C'$ and  $\ind C, \ind C' \subseteq \mathcal{M}$, we have 
	\begin{equation*}
	\ind \left(C \wedge C'\right) = \ind C \cap \ind C' \subseteq \mathcal{M}.
	\end{equation*}
	So $C \wedge C' \in \mathsf{S}(\mathcal{M})$.
	
	On the other hand, the join $C \vee C'$ consists of all subobjects of direct sums of objects in $C$ and $C'$. Thus, if $M \in \ind( C \vee C')$ then $M \in \sub \mathcal{M}$. By (\ref{EqSubinM}), $M \in \mathcal{M}$. So $C \vee C' \in \mathsf{S}(\mathcal{M})$ and $\mathsf{S}(\mathcal{M})$ is a sublattice of $\mathsf{S}(\mathcal{A})$.
\end{proof}

We get the following homomorphism between $\mathsf{S}(\mathcal{A})$ and a Cartesian product of sublattices of the form $\mathsf{(\mathcal{M})}$:

\begin{l1} \label{LemCartesianProduct}
	Let $\mathcal{A}$ be an abelian length category. Suppose that there is an index set $I$, and classes of indecomposable objects $\mathcal{M}_i$, $i \in I$ exist, so that 
	\begin{enumerate}[(M1)]
		\item 
		$\bigcup_{i\in I} \mathcal{M}_i = \ind A$
		
		\item $\mathcal{M}_i \cap \mathcal{M}_j = \emptyset$ for all $i, j \in I$ with $i \neq j$
		
		\item $\ind \sub \mathcal{M}_i = \mathcal{M}_i$ for all $i \in I$.
	\end{enumerate}
	Denote $\mathcal{M} = \{\mathcal{M}_i ~|~ i \in I\}$. Then
	\begin{equation*}
	f_{\mathcal{M}}: \mathsf{S(mod} ~A) \rightarrow \prod_{i\in I} \mathsf{S}(\mathcal{M}_i), \quad  C \rightarrow \prod_{i \in I} C_i
	\end{equation*}
	where $C_i$ is given by 
	\begin{equation*}
	\ind C_i = \ind C \cap \mathcal{M}_i
	\end{equation*}
	is a lattice homomorphism.
\end{l1}

\begin{proof}
	By Lemma \ref{LemSublattice} and (M3), $\mathsf{S}(\mathcal{M}_i)$ is a lattice for every $i \in I$ and the Cartesian product exists. We have to show that $f_{\mathcal{M}}$ preserves meets and joins. Take  $ C, C' \in \mathsf{S}(\mathcal{M}_i)$. 
	Then $f_{\mathcal{M}}$ preserves meets, since
	\begin{align*}
	\ind (C\wedge C')_i &= \ind (C \wedge C') \cap \mathcal{M}_i\\ &= (\ind C \cap \ind C') \cap \mathcal{M}_i\\&= \ind C_i \cap \ind C'_i\\ & = \ind (C_i \wedge C'_i),
	\end{align*}
	where the last equality holds by the definition of $\wedge$. Thus $(C \wedge C')_i = C_i \wedge C'_i$ and 
	\begin{equation*}
	f_{\mathcal{M}}(C \wedge C') = \prod_{i\in I} (C\wedge C')_i = \prod_{i\in I} (C_i \wedge C'_i)= \prod_{i \in I} C_i \wedge \prod_{i \in I} C'_i= f_{\mathcal{M}}(C) \wedge f_{\mathcal{M}}(C').
	\end{equation*}
	
	The function also preserves joins: For some object $X$, we have $X \in \ind (C \vee C')_i$ if and only if $X \in \mathcal{M}_i$ and there are objects $X_1, \dots, X_c \in \ind C \cup \ind C'$ for some $c, \in \mathbb{N}$ so that 
	\begin{equation*}
	X \subseteq \bigoplus_{k=1}^c X_k.
	\end{equation*}
	By Lemma \ref{LemSimpleObjects}, we can assume that for all $X_k$, $1 \le k \le c$, there is some simple $S \subset X$ with $S \subset X_k$. By (M1), there is some $j$ with $X_k \in \mathcal{M}_j$. We get $S \in \mathcal{M}_j$ and $S \in \mathcal{M}_i$ with (M3). Thus $X_1, \dots X_c \in \mathcal{M}_i$ and $X \in \ind (C_i \vee C'_i)$. So $\ind(C \vee C')_i \subset \ind(C_i \vee C'_i)$.
	
	To show the other direction, we suppose that $X \in \ind(C_i \vee C'_i)$. Then $X \in \mathcal{M}_i$ and $X \in C \vee C'$, since $C_i$ and $C'_i$ are subcategories of $C$ and $C'$ respectively. Thus, $X \in \ind(C \vee C')_i$.
	
	So $\ind (C_i \vee C'_i) = \ind (C \vee C')_i$ and $C_i \vee C'_i = (C \vee C')_i$. We get
	\begin{equation*}
	f_{\mathcal{M}}(C \vee C') = \prod_{i\in I} (C\vee C')_i = \prod_{i\in I} (C_i \vee C'_i)= \prod_{i \in I} C_i \vee \prod_{i \in I} C'_i= f_{\mathcal{M}}(C) \vee f_{\mathcal{M}}(C')
	\end{equation*}
	and $f_{\mathcal{M}}$ is a lattice homomorphism.
\end{proof}

Even better, $f_{\mathcal{M}}$ is an isomorphism:

\begin{p1} \label{PropCartesianProduct}
	Let $\mathcal{A}$ be an abelian length category and $\mathcal{M} = \{\mathcal{M}_i ~|~ i \in I\}$ be a family of classes of indecomposable objects that fulfil (M1) - (M3). Then $f_{\mathcal{M}}$ as defined in Lemma \ref{LemCartesianProduct}	is a lattice isomorphism between $\mathsf{S}(\mathcal{A})$ and $\prod_{i \in I} \mathsf{S}(\mathcal{M}_i)$.
\end{p1}

\begin{proof}
	By Lemma \ref{LemCartesianProduct}, $f_{\mathcal{M}}$ is a homomorphism between lattices. To show that $f_{\mathcal{M}}$ is an isomorphism, we need to prove that $f$ is injective and surjective.
	
	Suppose that $f_{\mathcal{M}}(C) = f_{\mathcal{M}}(C')$ for some $C, C' \in \mathsf{S(mod} ~A)$. Then 
	\begin{equation*}
	\prod_{i \in I} C_i = \prod_{i \in I} C'_i
	\end{equation*}
	and by (M2), we have $C_i = C'_i$ for all $i \in I$ . This means
	\begin{equation*}
	\ind C \cap \mathcal{M}_i = \ind C' \cap \mathcal{M}_i
	\end{equation*} 
	for all $i \in I$. By (M1), $\ind C = \ind C'$ and $f_{\mathcal{M}}$ is injective.
	
	Now take
	\begin{equation*}
	\prod_{i \in I} C_i \in \prod_{i\in I} \mathsf{S}(\mathcal{M}_i).
	\end{equation*} 
	Since all $C_i$ are subobject closed subcategories of $\mathcal{A}$, we have $C_i \in \mathsf{S} (\mathcal{A})$ for all $i \in I$. We will show that
	\begin{equation*}
	f_{\mathcal{M}}(\bigvee_{i \in I} C_i) = \prod_{i \in I} C_i.
	\end{equation*}
	It is obvious that $C_j \subseteq \left(\bigvee_{i\in I} C_i\right)_j $ for all $j \in I$ which implies
	\begin{equation*}
	\prod_{i \in I} C_i \subseteq f_{\mathcal{M}}(\bigvee_{i \in I} C_i).
	\end{equation*}	
	For the other direction, we need to show that  $\left(\bigvee_{i\in I} C_i\right)_j \subseteq C_j$ for all $j \in I$, which is equivalent to
	\begin{equation} \label{EqSurjective}
	\left(\ind \bigvee_{i \in I} C_i \right) \cap \mathcal{M}_j\subseteq \ind C_j. 
	\end{equation}
	Suppose that $X \in \left(\ind \bigvee_{i \in I} C_i \right) \cap \mathcal{M}_j$. Then there are objects $Y_i \in C_i$, so that 
	\begin{equation*}
	X \subseteq \bigoplus_{i \in I} Y_i.
	\end{equation*} 
	By Lemma \ref{LemSimpleObjects}, we can assume that for all $i \in I$, there is a simple object $S \subset X$ and $S \subset Y_i$.
	Using (M3), we get $S \in \mathcal{M}_j$ and $S \in \mathcal{M}_i$. So (M2) yields $I' = \{j\}$. Thus (\ref{EqSurjective}) holds, $f_{\mathcal{M}}$ is surjective and the proof is complete.
\end{proof}

\section{String algebras} \label{SecString}
	
	A special kind of quiver algebras are string algebras as described in \cite{BR}, Section 3:

\begin{d1}
	Suppose that $Q$ is a quiver and $I$ an ideal in $kQ$ which is generated by a set of zero relations. 
	
	Then $A= kQ/I$ is a \textit{string algebra} if and only if
	\begin{enumerate} \label{DefStringAlg}
		\item Any vertex of $Q$ is starting point of at most two arrows.
		
		\item Any vertex of $Q$ is end point of at most two arrows.
		
		\item Given an arrow $\beta$, there is at most one arrow $\gamma$ with
		$s(\beta) = e(\gamma)$ and $\beta \gamma \notin I$.
		
		\item Given an arrow $\gamma$, there is at most one arrow $\beta$ with
		$s(\beta) = e(\gamma)$ and $\beta \gamma \notin I$
		
		\item Given an arrow $\beta_1$, there is some bound $n(\beta_1)$ such that any path of the form $\beta_1 \beta_2 \dots \beta_{n(\beta_1)}$ contains a subpath in $I$.
		
		\item Given an arrow $\beta$, there is some bound $n'(\beta)$ such that any	path of the form $\beta_1 \beta_2 \dots \beta_{n'(\beta)}$ with $\beta_{n'(\beta)} = \beta$ contains a subpath in $I$.
	\end{enumerate}
\end{d1}

\begin{d1} \label{DefString}
	We can take the formal inverse $\beta^{-1}$ of an arrow $\beta$ by defining $e(\beta^{-1}):= s(\beta_n)$ and $s(\beta^{-1}):= e(\beta)$.
	
	A \textit{string} is a word $w = \beta_1 \beta_2 \dots \beta_n$ so that 
	\begin{itemize}
		\item $\beta_i$ is either an arrow or the inverse of an arrow for all $1 \le i \le n$
		\item $s(\beta_i) = e(\beta_{i+1})$ for all $1 \le i \le n$
		\item $w$ does not contain a relation in $I$
	\end{itemize}
	
	The multiplication of strings is analogous to the multiplication of paths of a quiver.
	
	A \textit{band} is a string $w = \beta_1 \beta_2 \dots \beta_n$ such that every power of $w$ is defined and does not contain a relation in $I$; furthermore $w$ may not be a power of a string $w' \neq w$.
\end{d1}

String algebras are especially useful, since their modules are well known, also from \cite{BR}, Section 3:

\begin{d1} \label{DefStringModule}
	Suppose that $w = \beta_1 \beta_2 \dots \beta_n$ is a string. Set $u(i) = e(\beta_{i+1})$, for 
	$0 \le i < n$, and $u(n) = s(\beta_n)$. 

	The \textit{string module} $M(w)$  is defined as the representation where for every $v \in Q_0$, the vector space $M(w)_v$ has as basis 
	\[
	\{z_i ~|~ u(i) = v\}
	\]
	with $z_i \neq z_j$ for $i \neq j$. 
	If $\beta_i$ is an arrow, then it defines the map $f_{\beta_i}(z_{i-1})=z_{i}$, otherwise $f_{\beta^{-1}_i}(z_i)=z_{i-1}$. For all other arrows $\alpha$, we have $f_{\alpha} = 0$.
	
	Now suppose that $w$ is even a band and $\phi: Z \rightarrow Z$ is an automorphism on a vector space over $k$. 

	The \textit{band module} $M(w, \phi)$ is defined as the representation with
	\[
	M(w, \phi)_v = \bigoplus_{e(\beta_{i+1})=v} Z_i
	\]
	where $Z_i = Z$.
	
	If $\beta_1$ is an arrow and $z \in Z_1$, then it defines the map $f_{\beta_1}(z) = \phi(z) \in Z_0$. If $\beta^{-1}_1$ is an arrow, then for $z \in Z_0$, $f_{\beta^{-1}_1}(z) = \phi^{-1}(z) \in Z_1$. 
	
	Let $2 \le i \le n$. If $\beta_i$ is an arrow and $z \in Z_i$, then $f_{\beta_i}(z) = z \in Z_{i-1}$; if $\beta^{-1}_i$ is an arrow and $z \in Z_{i-1}$, then $f_{\beta^{-1}_i}(z) = z \in Z_i$.
	
	For all other arrows $\alpha$, we have $f_{\alpha} = 0$.
\end{d1}

All modules over a string algebra are either string modules or band modules:

\begin{l1} \label{LemStringmodules}
	Let $A=kQ/I$ be a string algebra with a string $w = \beta_1 \beta_2 \dots \beta_n$.
	\begin{enumerate}[(a)]
		\item All $A$-modules are isomorphic to a string module or a band module
		
		\item Two string modules $M(w)$ and $M(w')$ are isomorphic if and only if $w=w'$ or $w' = w^{-1} := \beta^{-1}_n \beta^{-1}_{n-1} \dots \beta^{-1}_1$. 
		
		\item Two band modules $M(w, \phi)$ and $M(w', \phi')$ are isomorphic if and only if $\phi$ and $\phi'$ are similar and $w$ or $w^{-1}$ is a cyclic permutation of $w'$. 
		
		\item No band module is isomorphic to a string module.
		
	\end{enumerate}
\end{l1}

In the following section, we use a result from \cite{C-B}, p. 34 about morphisms between tree modules that reduces very nicely to monomorphisms between string modules: 

\begin{l1} \label{LemStringSubmodules}
	$M(w)$ is a submodule of $M(w')$ if and only if there are arrows $\alpha$, $\beta$ and strings $w_1, w_2$ so that $w'$ or $w'^{-1}$ is of the form
	\begin{equation*}
	w_1 \alpha^{-1} w \beta w_2
	\end{equation*}
	or
	\begin{equation*}
	w \beta w_2
	\end{equation*}
	or
	\begin{equation*}
	w_1 \alpha^{-1} w.
	\end{equation*}
\end{l1}

\section{The structure of the lattice} \label{SecStructure}

Let $A$ be Morita equivalent to $kQ/I$ for some quiver $Q$ and some admissible ideal $I$. This is always the case if $A$ is an algebra over a closed field $k$. If $A$ is of colocal type, then the lattice $\mathsf{S(mod} ~A)$ is relatively simple and can be described completely. 

This is actually a description of most distributive lattices of quiver algebras:

\begin{l1} \label{LemKronecker}
	Let $A$ be Morita equivalent to $kQ/I$ for some quiver $Q$ and some admissible ideal $I$. Suppose that $Q$ has the Kronecker quiver $\xymatrix{\bullet \ar@<0.5ex>[r] \ar@<-0.5ex>[r]& \bullet}$ as a subquiver. 
	
	Then $\mathsf{S}(\mod A)$ is not distributive.
\end{l1}

\begin{proof}
	There is a monomorphism
	\begin{equation*}
	Y :=\xymatrix{ k \ar@<0.5ex>^{\left[\begin{smallmatrix}
			1 \\ 0
			\end{smallmatrix}\right]}[r] \ar@<-0.5ex>_{\left[\begin{smallmatrix}
			0 \\ 1
			\end{smallmatrix}\right]}[r]& k^2} \rightarrowtail \xymatrix{k \ar@<0.5ex>^1[r] \ar@<-0.5ex>_0[r]& k} \oplus \xymatrix{k \ar@<0.5ex>^0[r] \ar@<-0.5ex>_1[r]& k} =: X_1 \oplus X_2.
	\end{equation*}
	Obviously, $Y \notin \sub X_1$ and $Y \notin \sub X_2$. Thus, $\mathsf{S}(\mod kQ)$ is not distributive by Proposition \ref{PropDistributive} and $\mathsf{S}(\mod A)$ is not distributive either.
\end{proof}

Together with Theorem \ref{ThShape2}, we get the following:

\begin{p1}
	Let $A$ be an Artin algebra which is Morita equivalent to $kQ/I$ for some quiver $Q$ and an admissible ideal $I$. 
\begin{enumerate}[(a)]
	\item The algebra $A$ is of colocal type if and only if $\mathsf{S(mod} ~A)$ is distributive and for every subquiver of $Q$ of the form
	\begin{equation*}
	\xymatrix{1 \ar^{\beta}[r] & 2 \ar@(dr, ur)_{\alpha}},
	\end{equation*}
	we have $\alpha \beta \in I$ or $\alpha^2 \in I$.
	 \item The algebra $A$ is of colocal type if and only if it is a string algebra and no vertex in $Q$ is starting point of more than one arrow.
\end{enumerate}
\end{p1}

\begin{proof}
First, we prove (a): By Proposition \ref{LemOneDirection}, $\mathsf{S(mod} ~A)$ is distributive if $A$ is of colocal type. The additional condition is fulfilled by Theorem \ref{ThShape2}, since it is implied by (C3).
	
On the other hand, suppose that $\mathsf{S(mod} ~A)$ is distributive and fulfils the condition above. For vertices $i, j$, let $S_i$, $S_j$ be the corresponding simple modules. Then $d^1_{S_i}(S_i, S_j)$ and $d^1_{S_j}(S_i, S_j)$ are both given by the number of arrows with starting point $i$ and end point $j$. By Lemma \ref{LemKronecker},
\begin{equation*}
d^1_{S_i}(S_i, S_j) = d^1_{S_j}(S_i, S_j) = 1.
\end{equation*}
 Thus Lemma \ref{LemTwoSimples} (b), \ref{LemFourSimples1} (b) and \ref{LemFourSimples2} (b) imply that $\mod A$ fulfils (C1) - (C3). By Theorem \ref{ThShape2}, $A$ is of colocal type.\\

To show (b), suppose that $A$ fulfils (C1) - (C3). This is equivalent to the following:
\begin{enumerate}
	\item  No vertex in $Q$ is starting point of more than one arrow.
	
	\item No vertex in $Q$ is end point of more than two arrows.
	
	\item Given an arrow $\beta$, there is at most one arrow $\gamma$ with
	$s(\beta) = e(\gamma)$ and $\beta \gamma \notin I$.
\end{enumerate}
Since $A$ is an Artin algebra, the quiver $Q$ must be finite. 

By Definition \ref{DefStringAlg}, it only remains to show that $I$ is an ideal generated by zero relations. Furthermore, every cycle in $Q$ is oriented, since every non-oriented cycle contains a vertex which is starting point of two arrows.

(1) also implies that the vertex $i$ belongs to a cycle, then there is exactly one arrow $\alpha$ with $s(\alpha) = i$ and it belongs to the cycle. So every connected component of $Q$ contains at most one cycle.

If for vertices $i$ and $j$, there is more than one path $\rho$ with $s(\rho) = i$ and $e(\rho) = j$, then all except for one of these paths contains an oriented cycle.

In fact, any relation which is not a zero relation is of the following form, where $\rho$ is an oriented cycle, $\rho'$, $\rho''$ are paths with $s(\rho') = e(\rho) = s(\rho) = e(\rho'')$, $a_1, \dots, a_n \in k$ and $\alpha_1 < \alpha_2 < \dots < \alpha_n \in \mathbb{N}$:
\begin{equation} \label{Relation}
a_1\rho'\rho^{\alpha_{1}}\rho'' + a_2 \rho'\rho^{\alpha_2}\rho'' + \dots + a_n\rho'\rho^{\alpha_n}\rho'' = 0.
\end{equation} 

Now, we use that $I$ is admissible: there must be some $t \in \mathbb{N}$, so that $\rho^t = 0$. So for every representation $V = (V_i, f_{\alpha})_{i \in Q_0, \alpha \in Q_1}$ of $Q$, there is some $m$ so that 
\begin{equation*}
0 = \Im f_{\rho^m} \subsetneq \Im f_{\rho^{m-1}} \subsetneq \dots \subsetneq  \Im f_{\rho}.
\end{equation*}
 We get $\rho^{\alpha_{1}} = \dots = \rho^{\alpha_n} = 0$, since otherwise
\begin{equation*}
\Im (f_{\rho'\rho^{\alpha_2}\rho''} + \dots + f_{\rho'\rho^{\alpha_n}\rho''}) \subseteq \Im f_{\rho'\rho^{\alpha_2}\rho''} \subsetneq \Im f_{\rho'\rho^{\alpha_1}\rho''},
\end{equation*}
which implies $f_{\rho'\rho^{\alpha_2}\rho''} + \dots + f_{\rho'\rho^{\alpha_n}\rho''} \neq f_{\rho'\rho^{\alpha_1}\rho''}$, a contradiction to (\ref{Relation}).
\end{proof}

Furthermore, we get some useful properties:

\begin{l1} \label{LemClassification}
	If $A = kQ/I$ for some quiver $Q = (Q_0, Q_1)$ with admissible ideal $I$ and $A$ is of colocal type, then the following holds:
	\begin{enumerate}[(a)]
		\item If $Q$ contains a cycle, this cycle is oriented.
		
		\item At most two paths are maximal under all paths without relations that end in $i$. If no or only one arrow ends in $i$, then there is only one such path.
		
		\item Every module in $\ind A$ is a string module.
		
		\item Every string is of the form
		\begin{equation} \label{EqW}
		w := \alpha^{-1}_{l_1} \alpha^{-1}_{l_1-1} \dots \alpha^{-1}_1 \beta_{1} \beta_{2} \dots \beta_{l_2},
		\end{equation}
		for some $l_1, l_2 \in \mathbb{N}_0$, and arrows $\alpha_1, \dots,  \alpha_{l_1}, \beta_1, \dots, \beta_{l_2}$ or of the form $e_m$ for some vertex $m$.
		
		\item Let $\underline{w}$ be defined as in (\ref{EqW}) and $\underline{w'}$ be a string. We have $M(w') \subseteq M(w)$ if and only if there are $1 \le j_1 \le l_1$ and $1 \le j_2 \le l_2$ so that $w'$ or $w'^{-1}$ is of the form
		\begin{equation*}
		\alpha^{-1}_{j_1} \alpha^{-1}_{j_1-1} \dots \alpha^{-1}_1 \beta_{1} \beta_{2} \dots \beta_{j_2},
		\end{equation*}
		or $w = e_m$ with $m = e(\alpha_1) = e(\beta_{1})$.
	\end{enumerate}
\end{l1}

\begin{proof}
	(a) Every non-oriented cycle contains a vertex which is starting point of two arrows.
	
	(b)Since $A$ is a string algebra, there are at most two arrows which end in $i$ by Definition \ref{DefStringAlg} (2). By \ref{DefStringAlg} (3), each of those arrows is part of only one maximal path that ends in $i$.
	
	(c) From definition \ref{DefString}, it is obvious that every band corresponds to a cycle without relations. Since $I$ is an admissible ideal, every oriented cycle of $Q$ contains a relation in $I$. By (a), $A = kQ/I$ has no band modules and $\ind A$ consists only of string modules.
	
	(d) There are no arrows $\alpha$, $\beta$ with $e(\beta^{-1}) = s(\beta) = s(\alpha)$. So no word contains a subword of the form $\alpha \beta^{-1}$.
	
	(e) This follows from Lemma \ref{LemStringSubmodules}.	
\end{proof}

We use Proposition \ref{PropCartesianProduct} to simplify the problem of describing $\mathsf{S(mod} ~A)$. First, we define a suitable family $\mathcal{M}$.

By Lemma \ref{LemClassification} (2), there are at most two maximal paths without relation that end in a vertex $m \in Q_0$:

\begin{d1}
	Suppose that there is at most one arrow $\alpha$ with $e(\alpha) = m$. Then there is only one path that is maximal under the paths without relation that ends in $m$. We denote its length with $k_m$ and set $l_m := 0$.
	
	If there are two arrows that end in $m$, there are two maximal paths. We denote their lengths with $k_m$ and $l_m$.
\end{d1}

\begin{d1}
	Let $A = kQ/I$ for some quiver $Q = (Q_0, Q_1)$, $m \in Q_0$, and $M_m := M(w_m)$ be the module with 
	\begin{equation*}
	w_m := \alpha^{-1}_{k_m} \alpha^{-1}_{l_1-1} \dots \alpha^{-1}_1 \beta_{1} \beta_{2} \dots \beta_{l_m}
	\end{equation*}
	so that $\alpha_{k_m} \dots \alpha_{1}$ and $\beta_{l_m} \dots \beta_{1}$ are the maximal paths that end in $m$.

	By Lemma \ref{LemStringmodules} (b) and Lemma \ref{LemClassification} (b), this module is well defined. 
	
	Furthermore, we define
	\begin{equation*}
	\mathcal{M}_m := \{ M \in \mod A ~|~ M \subset M_m \}.
	\end{equation*}
\end{d1}

\begin{l1} \label{LemMi}
	If $A$ is of colocal type, then
	\begin{equation*}
	\mathsf{S(mod} ~A) \cong \prod_{m\in Q_0} \mathsf{S}(\mathcal{M}_m).
	\end{equation*}
\end{l1}

\begin{proof}
	We need to prove that the sets $\mathcal{M}_m$, $m \in Q_0$ fulfil the conditions of Proposition \ref{PropCartesianProduct}:
	
	(M1) is fulfilled by Lemma \ref{LemClassification} (c), (d) and (e); (M2) and (M3) are fulfilled by Lemma \ref{LemClassification} (e).
\end{proof}

The lattices $\mathsf{S}(\mathcal{M}_m)$ for $m \in Q_0$ have a very simple description: they are all sublattices of Young's lattice, which is defined in \cite{Sta}, p. 58 and Example 3.4.4(b):

\begin{d1}
	Let $\lambda = (\lambda_1, \lambda_2, \lambda_3, \dots, \lambda_n)$ be  a partition of a natural number, ordered so that $\lambda_1 \ge \lambda_2 \ge \dots \ge \lambda_n$. The \textit{Young diagram} of $\lambda$ is an array of squares with $n$ rows and exactly $\lambda_i$ squares in the $i$-th row.
	
	These partitions form a lattice $Y$, ordered by the inclusion order on the Young diagrams. It is called \textit{Young's lattice}.
	
	Let $\lambda' := (\lambda'_1, \lambda'_2, \lambda'_3, \dots, \lambda'_n)$, suppose that $n \le n'$ and set $\lambda_i := 0$ for $i > n$. Then 
	\begin{equation*}
	\lambda' \wedge \lambda=\lambda \wedge \lambda' = (\min(\lambda_1, \lambda'_1), \dots, \min(\lambda_{n}, \lambda'_{n}) ))
	\end{equation*}
	and 
	\begin{equation*}
	\lambda' \vee \lambda=\lambda \vee \lambda' = (\max(\lambda_1, \lambda'_1), \dots, \max(\lambda_{n'}, \lambda'_{n'}) )).
	\end{equation*}
\end{d1}

\begin{e1}
	The Young diagram of the partition $(5, 3, 2, 1)$ has the following form:
	\begin{equation*}
	\ydiagram{5,3,2,1}
	\end{equation*}
\end{e1}

We will need the following lattices to describe $\mathsf{S}(\mathcal{M}_m)$ for $m \in Q_0$:

\begin{d1}
	Denote by $Y^{m,n}$ that sublattice of Young's lattice that contains exactly those partitions $\lambda = (\lambda_1, \lambda_2, \lambda_3, \dots, \lambda_{m'})$ where $m' \le m$ and $\lambda_i \le n$ for all $1 \le i \le m'$. Equivalently, we can define $Y^{m,n}$ as the lattice given by all Young diagrams with at most $m$ rows and at most $n$ columns.
\end{d1}

\begin{e1}
	The Hasse diagram of the lattice $Y^{3,3}$ is
	
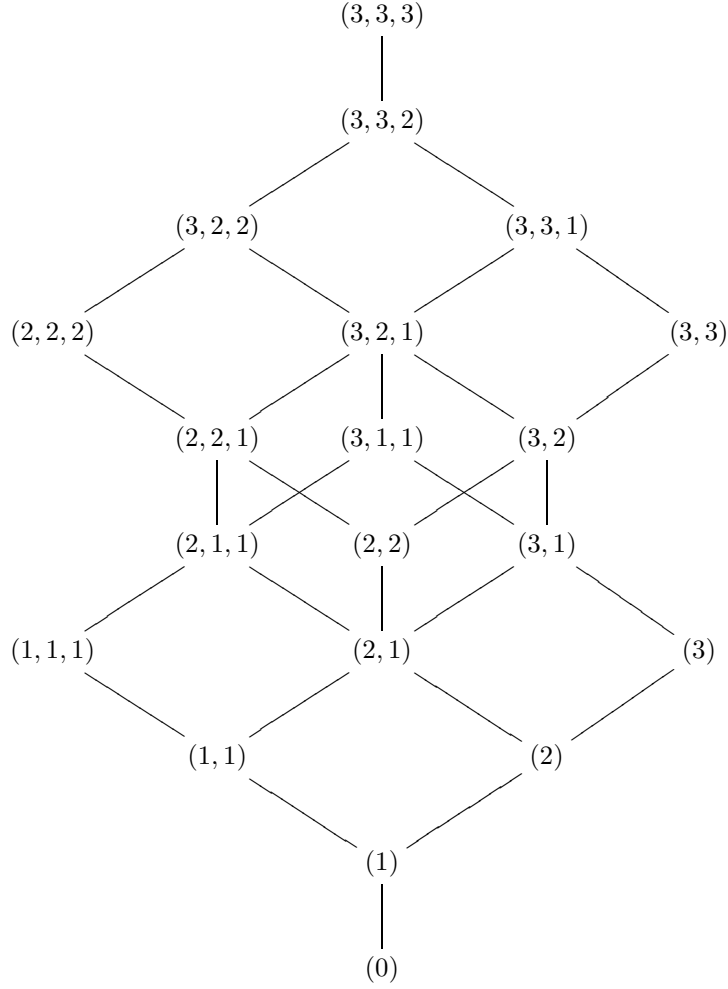
\begin{equation*}
	\xymatrix{
		&&(3,3,3)&&\\
		&& (3,3,2) \ar@{-}[u]&&\\
		& (3,2,2)\ar@{-}[ur]&& (3,3,1)\ar@{-}[ul]&\\
		(2,2,2)\ar@{-}[ur] & &(3,2,1) \ar@{-}[ul]\ar@{-}[ur]& &(3,3)\ar@{-}[ul]\\
		& (2,2,1)\ar@{-}[ul]\ar@{-}[ur]& (3,1,1)\ar@{-}[u]& (3,2)\ar@{-}[ul]\ar@{-}[ur]& \\
		&(2,1,1)\ar@{-}[ur]\ar@{-}[u]& (2,2) \ar@{-}[ul]\ar@{-}[ur]& (3,1)\ar@{-}[ul]\ar@{-}[u]&\\
		(1,1,1) \ar@{-}[ur]& &(2,1) \ar@{-}[ul]\ar@{-}[ur] \ar@{-}[u]& &(3)\ar@{-}[ul]\\
		& (1,1) \ar@{-}[ul]\ar@{-}[ur]& &(2) \ar@{-}[ul]\ar@{-}[ur]&\\
		& & (1) \ar@{-}[ul] \ar@{-}[ur]& &\\
		& & (0) \ar@{-}[u]& & & \\
	}
	\end{equation*}
\end{e1}

\begin{r1}
	Note that for $m, n \in \mathbb{N}$, we have $Y^{m, n} \cong Y^{n, m}$ and 
	\begin{equation*}
	Y^{1, n} \cong (\{0, 1, \dots, n\}, <) \cong Y^{n, 1}
	\end{equation*}
\end{r1}

Now, we can completely describe the distributive lattices $\mathsf{S(mod} ~A)$:  

\begin{t1}
	Suppose $A=kQ/I$ with quiver $Q = (Q_0, Q_1)$ and admissible ideal $I$. If $A$ is of colocal type, then 
	\begin{equation*}
	\mathsf{S(mod} ~A) \cong \prod_{m\in Q_0} Y^{k_m+1, l_m+1}.
	\end{equation*}
\end{t1}

\begin{proof}
	By Lemma \ref{LemMi}, 
	\begin{equation*}
	\mathsf{S(mod} ~A) \cong \prod_{m\in Q_0} \mathsf{S}(\mathcal{M}_m).
	\end{equation*}
	If only one path $\alpha_1 \alpha_{2} \dots \alpha_{k_m}$ ends in $m$, it is obvious from Lemma \ref{LemClassification} (e) that we can order the modules in $\mathcal{M}_m$ the following way:
	\begin{equation*}
	M(e_m) \subseteq M(\alpha^{-1}_1) \subseteq M(\alpha^{-1}_{2} \alpha^{-1}_1) \subseteq \dots \subseteq M(\alpha^{-1}_{k_m} \alpha^{-1}_{k_m-1} \dots \alpha^{-1}_1)
	\end{equation*}
	Thus
	\begin{equation*}
	\mathsf{S}(\mathcal{M}_m) \cong  (\{0, \dots, k_m+1\}, <)= Y^{k_m+1, 1}.
	\end{equation*}
	If there are two maximal paths without relations $\alpha_{1} \alpha_{2} \dots \alpha_{l_2}$ and $\beta_{1} \beta_{2} \dots \beta_{l_2}$ that end in $m$, then by \ref{LemClassification} (e) all modules in $\mathcal{M}_m$ are of the form $M(0,0) := M(e_m)$ or $M(i,j) := M(w_{i j})$ with
	\begin{equation*}
	w_{ij} := \alpha^{-1}_{i} \alpha^{-1}_{i-1} \dots \alpha^{-1}_1 \beta_{1} \beta_{2} \dots \beta_{j},
	\end{equation*}
	and $0 \le i \le k_m$, $0 \le j \le l_m$ and $1 \le i+j$. Furthermore, $M({i,j}) \subseteq M({i',j'})$ if and only if $i \le i'$ and $j \le j'$.
	
	So for a submodule closed subcategory $C \in \mathsf{S}(\mathcal{M}_m)$, there are some $j_0, j_2, \dots, j_{\alpha}$ so that $M({i, j}) \in C$ if and only if $M({i, j}) \subset M({i, j_i})$, which is equivalent to $j \le j_i$. In particular, $j_{i} \ge j_{i+1}$ for $0 \le i < \alpha$, since $M({i,j_{i+1}}) \subset M({i+1, j_{i+1}})$.

	Because all modules in $C$ are submodules of $M({k_ml_m})$, we get
	\begin{equation*}
		 l_m \ge j_0 \ge j_1 \ge \dots \ge j_{\alpha} \ge 0
	\end{equation*}
	and $\alpha \le k_m$. 

	We define
	\begin{equation*}
	\lambda_C := (j_0+1, j_1+1, \dots j_{\alpha}+1).
	\end{equation*}
	Then $\lambda_C$ is well-defined and
	\begin{equation*}
	f: \mathsf{S}(\mathcal{M}_m) \rightarrow Y^{k_m+1, l_m+1}, C \rightarrow \lambda_C
	\end{equation*}
	is obviously injective and surjective. We need to prove that $f$
	is a lattice homomorphism, that is, that it preserves joins and meets: Since $\mathsf{S}(\mathcal{M}_m)$ is distributive, for any two categories $C_1, C_2 \in \mathsf{S}(\mathcal{M}_m)$, we have
	\begin{equation*}
	\ind (C_1 \wedge C_2) = \ind C_1 \cap \ind C_2
	\end{equation*}
	and 
	\begin{equation*}
	\ind (C_1 \vee C_2 ) = \ind C_1 \cup \ind C_2
	\end{equation*}
	by Proposition \ref{PropDistributive}.
	From the definition of the joins and meets in $Y^{k_m+1, l_m+1}$, it is clear that $f$ preserves them.
\end{proof}

\section*{Acknowledgements}

I would like to thank Claus Ringel for his numerous comments and his advice.


\begin{thebibliography}{99}
		
		\bibitem{ASS} I. Assem, D. Simson, A. Skowro\'nski, \textit{Elements of the Representation Theory of Associative Algebras. Volume 1}, Cambridge University Press, Cambridge 2006.

		\bibitem{Ben} D. J. Benson, \textit{Representations and Cohomology I: Basic Representation Theory of Finite Groups and Associative Algebras}, first paperback edition, Cambridge University Press, Cambridge 1998.
		
		\bibitem{BR} M C. R. Butler, C. M. Ringel, \textit{Auslander-Reiten sequences with few middle terms and applications to string algebras}, Communications in Algebra 15(1\&2) (1987), 145-179.
	
		\bibitem{C-B} W. Crawley-Boevey, \textit{Infinite-dimensional modules in the representation theory of finite-dimensional algebras}, Canadian Mathematical Society, Conference Proceedings, 23 (1998), 29-54.
		
		\bibitem{Gab} P. Gabriel, \textit{Indecomposable representations II}, Symposia Mathematica 11 (1973), 81-104.
		
		
		\bibitem{Got} A. Gottwald, \textit{Cofinite Submodule Closed Categories and the Weyl group}, arXiv:1705.04280
		
		\bibitem{IK} S. B. Iyengar, H. Krause, \textit{The Bousfield lattice of a triangulated category and stratification}, Mathematische Zeitschrift, Vol. 273, Issue 3-4 (2013), 1215-1241.

		\bibitem{KP} H. Krause, M. Prest \textit{The Gabriel-Roiter Filtration of the Ziegler Spectrum}, The Quarterly Journal of Mathematics, Vol. 64 Issue 3 (2013), 891-901.
		
		\bibitem{ORT} S. Oppermann, I. Reiten, H. Thomas, \textit{Quotient closed subcategories of quiver representations}, Compositio Mathematica, 151 (2015), 568-602.
		
		\bibitem{Osb} M. S. Osborne, \textit{Basic Homological Algebra}, Springer Verlag, New York 2000.
		
		\bibitem{Rin} C. M. Ringel, \textit{Minimal Infinite Submodule-closed Subcategories}, Bulletin des Sciences Math\'ematiques 136 (2012), 820-830.

		\bibitem{Ste} B. Stenstr\"om, \textit{Rings of Quotients}, Springer Verlag, Berlin, Heidelberg 1975.
		
		\bibitem{Sta}  R.P. Stanley, \textit{Enumerative Combinatorics, Volume 1}, second Cambridge edition, Cambridge University Press, New York 2012
		
		\bibitem{Sum} T. Sumioka, \textit{Tachikawa's Theorem on Algebras of Left Colocal Type}, Osaka J. Math. 21 (1984), 629-648
		
		\bibitem{Tac} H. Tachikawa, \textit{On Rings for which Every Indecomposable Right Module has a Unique Maximal Submodule}, Math. Zeitschr. 71 (1959), 200-222.
		
		\bibitem{Zim} A. Zimmermann, \textit{Representation Theory. A Homological Algebra Point of View},  Springer International Publishing Switzerland 2014.
	\end{thebibliography}
\end{document}